\documentclass[11pt]{article}
\usepackage{graphicx}
\usepackage{amsfonts}
\usepackage{amssymb}
\usepackage{amsmath,amsthm}
\usepackage{a4wide}

\usepackage{xcolor}

\newcommand{\dis}{\displaystyle}

 \newcommand{\nn}{\nonumber}

\numberwithin{equation}{section}
\theoremstyle{plain}
\newtheorem{theorem}{Theorem}[section]
\newtheorem{lemma}[theorem]{Lemma}
\newtheorem{corollary}[theorem]{Corollary}
\newtheorem{proposition}[theorem]{Proposition}
\theoremstyle{definition}
\newtheorem{definition}{Definition}[section]
\newtheorem{remark}{Remark}[section]

\def\si{{\sigma}}
\def\la{{\lambda}}

\def\ga{{\gamma}}
\def\Ga{{\Gamma}}

\def\N{{\mathbb N}}

\def\eps{{\epsilon}}

\def\square{\ifmmode\sqr\else{$\sqr$}\fi}
\def\sqr{\vcenter{
         \hrule height.1mm
         \hbox{\vrule width.1mm height2.2mm\kern2.18mm\vrule width.1mm}
         \hrule height.1mm}}                  

\newcommand{\col}[1]{\textcolor[rgb]{0,0,0}{#1}}
\newcommand{\colo}[1]{\textcolor[rgb]{0,0,0}{#1}}
\newcommand{\cola}[1]{\textcolor[rgb]{0,0,0}{#1}}

 \title{Global solutions of a free boundary problem\\
 via mass transport inequalities }
\author{
Gioia Carinci$^{\textup{{\tiny(a)}}}$,
Anna De Masi $^{\textup{{\tiny(b)}}}$,\\
Cristian Giardin{\`a} $^{\textup{{\tiny(a)}}}$,
Errico Presutti $^{\textup{{\tiny(c)}}}$.\\\\
{\small $^{\textup{(a)}}$ Dipartimento di Dipartimento di Scienze Fisiche, Informatiche e Matematiche,}\\
{\small Universit\`a di Modena e Reggio Emilia, via Campi 213/b, 41125 Modena, Italy}
\\
{\small $^{\textup{(b)}}$ Dipartimento di Ingegneria e Scienze dell'Informazione e Matematica,}\\
{\small Universit\`a di L'Aquila, via Vetoio 1, 67100 L'Aquila, Italy}\\
{\small $^{\textup{(c)}}$ GSSI, viale F. Crispi 7, 67100 L'Aquila, Italy}\\
{\small }\\
}

\date{\today}

\begin{document}

%
\maketitle

\begin{abstract}


\col{
We study a free boundary problem
which arises
as the continuum version of a stochastic particles system in the context
of Fourier law.  Local existence and uniqueness of the classical solution
are well known in the literature of {free boundary}
 problems.
We introduce the notion of  generalized solutions
(which extends that of classical solutions when the latter exist) and prove
global existence and uniqueness of  generalized solutions for a large class of
initial data.  The proof is obtained by  characterizing a generalized solution
as the unique element which separates suitably defined lower and upper barriers
in the sense of  mass transport inequalities.
}

\end{abstract}

%
%
%
%
%
%
%
%
%
%
%



\vskip1cm

\section{Introduction and results}
\label{sec:1}

When we heat a metal bar from one side and cool it from the other
we produce a heat flow whose law  has been established by J.B.J. Fourier
\cite{fourier}.  We shall consider here a variant where
there is a mass  rather than a heat flow  across
the system.
The analogue of the Fourier's law of heat
conduction is then the Fick's law of mass transport,
which is also described
by the heat equation.  Since the transversal direction to the flow is not
relevant we model our system as one dimensional.
A classical experiment of mass transport
is the following: the system is confined in the interval $[0,X]$,
$X>0$, we inject mass into the system
from the left boundary 0 at rate $j>0$, the mass diffuses in
$[0,X]$ and it is removed once it reaches the right boundary $X$.
Supposing a constant conductivity (set  equal to 1), then, according to the Fick's law,
the mass density $\rho(r,t)$   solves the equation
  \begin{equation}
    \label{0.1}
 \frac{\partial \rho}{\partial t} = \frac 12
 \frac{\partial^2 \rho}{\partial r^2}
   \end{equation}
in $[0,X]$ with initial condition $\rho(r,0)=\rho_0(r)$ and boundary conditions
  \begin{equation}
    \label{0.2}
-\frac 12 \frac{\partial \rho}{\partial r}(0,t) = j,\quad
  \rho(X,t)= 0 
   \end{equation}

The Fick's law described by \eqref{0.1} and  \eqref{0.2} then predicts how much
mass we are removing from the system at the right boundary: the removal rate is
  \begin{equation}
    \label{0.3}
-\frac 12 \frac{\partial \rho}{\partial r}(X,t)\;.
   \end{equation}
We study in this paper the case when the endpoint $X$ is no longer fixed but
free to slide on the line pushed
by the  mass flow.  We exploit this to ensure that the total mass is conserved.
Thus $X_t$ is the position
of a moving reservoir, whose spatial evolution guarantees that at each time
there is a perfect balance between the amount of mass per unit time  that is injected
into the system, equal to $j$, and the amount that is extracted, which should therefore be
also equal to $j$.
%
%
%
%
%
Fick's law in this case is then  again \eqref{0.1} but in a time changing interval
$[0,X_t]$ with initial condition  $\rho(r,0)=\rho_0(r)$ and boundary conditions
  \begin{equation}
    \label{0.4}
-\frac 12 \frac{\partial \rho}{\partial r}(0,t) = j,\;
-\frac 12 \frac{\partial \rho}{\partial r}(X_t,t) = j,\quad
  \rho(X_t,t)=0\;.
   \end{equation}
This is a free boundary problem  (FBP) where the unknown are both the solution
of the equation and the interval where the equation is set.

There exist particle versions of the system, in particular the paper
\cite{CDGP}  which has actually motivated the present work.  There
are  also several works on particle systems which model free boundary
problems, we quote here a few of them which are closer to our system,
\cite{CDGP,CDGP2, dfp, dptv1,dptv2,dpt,GQ,lacoin}.

Free boundary problems  for the linear heat equation have been
thoroughly  investigated, a celebrated example is the Stefan problem,
i.e. the description of the interface between ice and water
in the presence of a first-order phase transition \cite{stefan}.
As we shall see the space derivative of the solution
of \eqref{0.1}--\eqref{0.4}
is related to the  classical Stefan problem.
In one dimension the theory \col{of Stefan FBP} is
very rich and detailed results are available
\col{\cite{douglas, fasano, friedman, luckhaus, visintin}}.

%
%
%
In particular, \col{under appropriate assumptions on the
initial datum,
global (in time) existence and uniqueness theorems are known.
However in the general case only existence results local in time are available due to the
appearance of singularities for the classical solutions.}

The aim of this paper is to study the extension of the \col{FBP} solution
past the singularities. \col{By a combination of probabilistic methods and results
in the context of mass transport theory we provide {\em global} existence
and uniqueness. The result for all times is achieved by introducing a new notion of
{\em generalized} solution obtained as the weak limit of the classical solutions of
an auxiliary FBP problem with a relaxed mass constraint. The scheme is applicable to
a large set of initial conditions, including bounded and integrable initial data.}
%
%

\subsection{The free boundary problem.}

{In Definition \ref{defin1.1} }below 
we give the classical formulation of the problem.

  \begin{definition}[The FBP \ref{defin1.1}]
    \label{defin1.1}

The pair  $(X_t, \rho(\cdot ,t))$ is a classical solution of the FBP  \ref{defin1.1}
in the time interval $[0,T)$ with initial datum $(X_0,\rho_0(\cdot))$ 
if it satisfies:
i) $X_t\in C^1([0,T), \mathbb R_+)$ is strictly positive and $X_{t=0}=X_0$;
ii) for each $t\in [0,T)$, $\rho(\cdot,t)\in
C^2((0,X_t), \mathbb R_+)$ and it has limits with its derivatives at $0$ and $X_t$; moreover for each $r\in [0,X_t]$, $\rho(r,t)$
is differentiable in $t$; 
iii) the following equations \cola{(with $j$ a positive parameter)} are pointwise satified
\begin{eqnarray}
\label{2.4a}
\frac{\partial \rho}{\partial t}(r,t)  &  = & \frac 12 \frac{\partial^2 \rho}{\partial r^2}(r,t)
\qquad\qquad t\in [0,T), r\in (0,X_t)\\ 
\label{2.4aa}
 \rho(r,0) &=&\rho_0(r),{\qquad\qquad\qquad r\in(0,X_0)}\\
\label{2.4b}
\frac{\partial \rho}{\partial r}(0,t)& = & -2 j
\qquad\qquad\qquad\quad t\in [0,T) \\
\label{2.4c}
\rho( X_t,t) & = & 0
\qquad\qquad\qquad\qquad\; t\in [0,T) \\
\label{2.4d}
\frac{\partial \rho}{\partial r}(X_t,t) & = & -2j
\qquad\qquad\qquad\quad t\in [0,T)
\end{eqnarray}
By \eqref{2.4aa} the initial data  $(X_0,\rho_0(\cdot))$ of classical solutions must be as regular
as above and satisfy \eqref{2.4b}--\eqref{2.4d} with $t=0$.
%
\end{definition}

\smallskip

\begin{remark}
As already mentioned in the beginning of the paper
$\rho(\cdot,t)$
is a non-negative function
that we interpret as a mass density which is
concentrated in the time-varying interval $[0,X_t]$.   $ X_t $ is
then
the right edge of the mass distribution and it is called free boundary
because it is itself part of the problem.

The condition \eqref{2.4b}
describes a constant incoming
mass flow (represented by $-\frac 12\frac{\partial \rho}{\partial r}(0,t)$) through
the origin, \col{with mass injected}  at rate $j >0$,
and the condition \eqref{2.4d} states that
there is an outgoing mass flow, $-\frac 12 \frac{\partial \rho}{\partial r}(X_t,t)$, through $X_t$
which is also equal to $j$.  Inside $(0,X_t)$
the mass diffuses freely according to the linear heat equation \eqref{2.4a},
thus
in a classical solution
the total mass is conserved
\begin{equation}
\label{1.4d}
 \int_{0}^{X_t}\rho(r,t)\,dr  =   \int_{0}^{X_0}\rho(r,0)\,dr\qquad\quad\quad\;\;\; t\in (0,T)\;.
\end{equation}
\col{This} can be formally seen by
differentiating the left hand side and using  \eqref{2.4a}---\eqref{2.4d};
\col{moreover one can check that classical solutions of FBP \ref{defin1.1}
coincide with classical solutions of  the FBP problem that is obtained  replacing
condition \eqref{2.4d} by \eqref{1.4d}}.
\end{remark}

\begin{remark}
According to the Fick's law once we fix the density values at the boundaries
we find that the stationary state is  a  linear function.  This suggests
that in our case the stationary
solutions are
$$
(X_t,\rho(r,t)) := \left(\frac{a}{2j},\col{\rho_{a}}(r)\right)
$$
where $a>0$ and
\begin{equation}
\label{1.6}
\col{\rho_{a}}(r) = a-2j r,\qquad 0\le r\le \frac{a}{2j}\;.
\end{equation}
We do not have a unique stationary state but rather the whole family
$\{\rho_{a}, a>0\}$,
this is because we do not fix the density at the endpoints but only the current:
physically this means that we deal with a ``current reservoir'' while usually
the Fick's law is applied to cases where there are
reservoirs which fix the densities at the endpoints.

We shall later use  the functions $\col{\rho_{a}}(r)$
to bound the solutions of the problem \eqref{2.4a}---\eqref{2.4d}  via
mass transport inequalities.
\end{remark}

\vskip.5cm

Local existence of classical solutions of  FBP \ref{defin1.1}
follows from the  literature \col{of Stefan problems}:

\medskip

\begin{theorem}[Local existence]
\label{localexist}
Suppose the initial datum $(X_0,\rho_0)$ \cola{satisfies}: $X_0>0$,
$\rho_0  \in C^3([0,\col{X_0}), \mathbb R_+)$, $\rho_0$
has limit with its derivatives at $\col{X_0}$ and $\rho_0(X_0)=0$, $\frac{d\rho_0}{dr}(X_0) = -2j$,
\col{$\frac{d\rho_0}{dr}(0) = -2j$}.
Then there \col{exists} $T>0$ and a classical solution
$(X_t,\rho(\cdot,t))$, $t\in [0,T)$, with initial datum $(X_0,\rho_0)$.

\end{theorem}

\medskip

\noindent
{\bf Proof.} Consider the classical Stefan problem in $ t\in [0,T), r\in (0,X_t)$:
 \begin{eqnarray}
\label{2.2}
&&\frac{\partial v}{\partial t}   =
\frac 12 \frac{\partial^2 v}{\partial r^2} \;\; 
\quad  v(r,t)\Big|_{r=0,X_t} = 0,\;\; v(r,0) =- \frac 12 \frac{\partial \rho_0(r)}{\partial r}  - j \nn
 \\&& \frac{dX_t}{dt} = -(2j)^{-1} \frac{\partial v(r,t)}{\partial r}\Big|_{r=X_t}
\end{eqnarray}
which is formally obtained from  \eqref{2.4a}---\eqref{2.4d} by setting
$\dis{v(r,t):=- \frac 12 \frac{\partial \rho}{\partial r}}(r,t)-j$;
the equation for $X_t$ being obtained by differentiating the identity
$\rho(X_t,t)=0$.   Local existence for \eqref{2.2} is proved in
\cite{fasano}--\cite {fp4}.

Given $X_t$ and  $v(r,t)$ satisfying \eqref{2.2} we set
 \begin{eqnarray}
\label{2.3}
&&\rho(r,t) = 2\int_{r}^{X_t} \Big(v(r',t)+j\Big)dr'
\end{eqnarray}
and check that \eqref{2.4a}---\eqref{2.4d} are satisfied by $(X_t,\rho(\cdot,t))$.
Non negativity of $\rho(\cdot,t)$ follows from the maximum principle, see also Section \ref{sec:4} where we express $\rho(\cdot,t)$ in terms of Green functions. \qed

\medskip

\col{Uniqueness of the local classical solution for the  Stefan problem \eqref{2.2} is also known in the literature,
however this does not immediately imply uniqueness  for
FBP \ref{defin1.1} which will be proved in Theorem \ref{exist}.}

\medskip

Following Fasano and Primicerio \col{(see e.g. \cite{fasano})} we say that  if $v(r,0) \ge 0$
then \eqref{2.2} has a ``sign specification''.  With a sign specification
the solution
is global while if
there is no sign specification in general we only have local existence with
examples where singularities do appear.  The analysis of their structure  is a very interesting
and much studied problem, see for instance \cite{Crank}, \cite{fp5}, \cite{fp6},
\cite{Ockendon}.

Our main goal in this paper is the extension of the solution
past the singularities. We thus want to introduce a notion of ``generalized solutions''
which reduces to classical solutions in the smooth case
and discuss global existence and uniqueness
of such generalized solutions.

\vskip1cm

\subsection{Quasi-solutions and generalized solutions.}

We define a generalized solution of the FBP \ref{defin1.1} as the weak limit of solutions of
approximate problems  calling the latter ``quasi-solutions''.

\begin{definition}[Quasi-solutions]
\label{defin1.2}

Let $\rho_0\in L^{\infty}(\mathbb R_+,\mathbb R_+)\cap L^1(\mathbb R_+,\mathbb R_+)$
be
such that:
 \begin{equation}
\label{1.7}
0<R(\rho_0) : = \inf\Big\{ r: \int_r^\infty \rho_0(r')\,dr' =0\Big\} < \infty
\end{equation}
Then $(X_t, u(\cdot,t),\eps)$,
  $t\in [0,T)$, $T>0$,
  $\eps>0$, is
a quasi-solution of the FBP \ref{defin1.1}
in the time interval $[0,T)$ with initial datum $\rho_0$ and accuracy parameter $\eps$ if
the following conditions are satisfied:

\begin{itemize}
\item  $X_t$ is strictly positive, Lipschitz continuous and piecewise $C^1$ (with finitely many
discontinuities of the derivative)

\item
$u(r,t)$, $t\in [0,T)$, $r\in [0,X_t]$,
is smooth (in the sense of (ii) of Definition \ref{defin1.1})
and it solves \eqref{2.4a}, \eqref{2.4b}, \eqref{2.4c} for all $t\in (0,T)$.

\item
The condition \eqref{2.4aa} is replaced by $\dis{\int | u(r,0) - \rho_0(r)| dr \le \eps}$.

\item
The condition \eqref{1.4d} is replaced by
\begin{equation}
\label{1.9}
\sup_{t\le T} \Big| \int_0^{X_t}u(r,t)\,dr -  \int_0^{X_0}u(r,0)\,dr\Big| \le \eps\;.
\end{equation}

\end{itemize}

\noindent
A  sequence $\{(X^{(n)}_t,\rho^{(n)}(\cdot,t),\eps_n)$, $t\in [0,T)\}$  of quasi-solutions in $[0,T)$
with initial datum $\rho_0$ is called ``optimal'' if $\dis{\lim_{n\to \infty} \eps_n = 0}$.

\end{definition}

\begin{definition}[Generalized solutions]
\label{defin1.2bis}
Let {$\rho_0$ be} as in Definition \ref{defin1.2} and $T>0$. Then
$\rho(r,t)$, $r\in \mathbb R_+$, $t\in [0,T)$, is a generalized solution in $[0,T)$ with initial datum $\rho_0$
of the FBP \ref{defin1.1}
if there exists an  {optimal}  sequence {$(X^{(n)}_t,\rho^{(n)}(\cdot,t),\eps_n)$, $t\in [0,T]$,
of quasi-solutions in \cola{the time interval} $[0,T)$
with initial datum $\rho_0$ such that for all $t\in [0,T)$
   \[
 \lim_{n\to \infty} \ \rho^{(n)} (\cdot,t) = \rho(\cdot,t) \quad \mbox{weakly}
   \]}

\end{definition}

\medskip
\noindent
The main result in this paper is a proof of {\em global existence and uniqueness of generalized solutions with a variational
representation of the generalized solution as the unique separating element of lower and upper barriers.}
{This will be explained in Section \ref{sec:1.3} below}.

%
%
%
\medskip

{
\begin{theorem}[Existence and uniqueness]
\label{exist}
For any $\rho_0$ as in Definition \ref{defin1.2} and  any  $T>0$ {the following holds}.
\begin{itemize}
\item[ (a)] There exists an optimal sequence   {
of quasi-solutions
 in $[0,T)$ with  initial datum $\rho_0$.

\item[(b)] Any optimal sequence
of quasi-solutions in $[0,T)$  with  initial datum $\rho_0$
converges weakly 
to a limit which
is (by definition) a generalized solution   of the FBP \ref{defin1.1}.}

\item[(c)]
There exists a function
$\bar \rho=\{\bar \rho(r,t)$, $r\in \mathbb R_+$,
$t\in \mathbb R_+\}$, continuous in $(r,t)$,
such that if $u(\cdot,t)$ is a
generalized solution of the FBP \ref{defin1.1} in  $[0,T)$ with initial datum $\rho_0$ then $u(\cdot,t)=\bar \rho(\cdot,t)$ for all $t\in[0,T)$.


\end{itemize}
\end{theorem}
}
\vskip.5cm

\begin{remark}
From  Theorem \ref{exist} we thus conclude that all optimal sequences of quasi-solutions
(in any  interval $[0,T)$ and with initial datum $\rho_0$) converge to the function $\bar \rho$, which is then the unique generalized solution.  Moreover since any classical solution of the FBP \ref{defin1.1} is also an optimal sequence of quasi-solutions with $\eps_n\equiv 0$ in the time interval where it exists then  it is  also a generalized solution. Thus when the unique generalized solution is smooth (in the sense of (ii) of Definition \ref{defin1.1})
it is a classical solution.

Finally observe that the FBP \ref{defin1.1} can be simulated by
the stochastic particle evolution studied in \cite{CDGP} which provides a discrete
approximation convergent in the macroscopic hydrodynamic limit.
\end{remark}

\vskip.3cm

{In the next section we give a characterization of the function $\bar \rho$ as the unique separating element between barriers. }

\vskip.5cm

\vskip1cm

\subsection{A variational representation of the evolution.}
\label{sec:1.3}

\col{The key} to the proof of Theorem \ref{exist}
is  variational as it exploits the
monotonicity properties (in the sense of mass transport) of the FBP \ref{defin1.1}
by introducing lower and upper barriers.

The notion of barriers for the construction of solutions
of partial differential equations
is well known   \cite{L, friedman2} (see also
\cite{DG} in the context of motion by mean curvature).
In our case upper and lower barriers will provide upper and lower
bounds with respect to an order based on mass transport.

We construct the barriers by quantizing the way to add
and remove mass to the system: we add
quanta of mass ($=j\delta$) at the origin at discrete times,
say $k\delta$, $\delta>0$, $k\in \mathbb N$,
and simultaneously we
remove the same amount of
mass at a left neighbor of the edge. Since at each step a finite mass is added all
at the origin we shall need to enlarge the functional space we work with
to include the Dirac delta (cfr. Definition \ref{def-u}).
In the time intervals  $(k\delta,(k+1)\delta)$ we
evolve using the linear heat equation on $\mathbb{R}_{+}$
with Neumann condition at 0.
We will obtain the upper barrier if
the addition/removal starts at time 0, and  the lower
barrier  if it starts at time $\delta$.
We shall prove,  see Theorem \ref{thmineq}, that
there is a unique element which separates the upper
and lower barriers constructed in this way.
The discrete
scheme by which the barriers are defined was previously used in \cite{CDGP}
and introduced in \cite{dfp} to study a similar particle model.

We define the barriers  in Definition \ref{defin:barriers} below after some preliminaries.

\medskip

\begin{definition}
\label{def-u}
(The space $\mathcal U$). $\mathcal U$ is the space of positive Borel measures on $\mathbb R_+$
which are sum of $cD_0$,  $c\ge 0$, $D_0$ 
the Dirac delta at 0, and $\rho\,dr$, with $\rho\in L^\infty(\mathbb R_+,\mathbb R_+)
\cap L^1(\mathbb R_+,\mathbb R_+)$.
We shall denote its elements by $u$ and by an abuse of notation, write
$u = c_uD_0+\rho_u$.  For $u$ and $v$ in $\mathcal U$ we call
  \begin{equation}
\label{3.2a}
|u-v| = |c_u-c_v| D_0 +|\rho_u-\rho_v|,\quad |u-v|_1= |c_u-c_v|   +\int |\rho_u-\rho_v|\,dr
   \end{equation}
We further define
  \begin{equation}
\label{3.2}
F(r;u)=  c_u\mathbf 1_{\{0\}}(r)+\int_r^\infty \rho_u(r')\,dr'
   \end{equation}
where
$\mathbf
1_{\{0\}}(r)=0$ unless $r=0$, in which case it is equal to 1. Note that $F(0;u)$ is the total mass of the measure $u$
and $F(0;|u-v|)=|u-v|_1$ is the total variation norm of $u-v$.
We call $\mathcal U_\delta$, $\delta \ge 0$, the following subset of $\mathcal U$:
%
  \begin{equation}
\label{3.3}
\mathcal U_\delta:=\Big\{u\in \mathcal U:\; F(0;\rho_u)> j\delta\Big\}
   \end{equation}

   \end{definition}

   \vskip.5cm
\begin{remark}
\noindent Observe that  $F(r;u)$ is a non increasing
function of $r$ which starts at 0 from
$F(0;u)$ which is the total mass of $u$.
If $u$ has compact support  $F(r;u)=0$ definitively
\col{and, in agreement with \eqref{1.7}, we can define the ``edge''
$R(u)$ as}
   \begin{equation}
\label{4.2a}
R(u):= \inf\{r: F(r;u)=0\}
   \end{equation}
\end{remark}

\medskip

\begin{definition}
 (The cut and paste operator).
The {\em cut-and-paste} operator
$K^{(\delta)}: \mathcal U_\delta\to \mathcal U$ is defined as
   \begin{equation}
\label{3.4}
K^{(\delta)} u = j\delta D_0+
\mathbf 1_{[0,R_\delta(u)]}u,\qquad R_\delta(u):= \inf\{r:F(r;u)= j\delta\}
   \end{equation}
\end{definition}

\medskip

\begin{definition}
(The free evolution).
Call $G_t^{\rm neum}(r,r')$, $r,r' \in \mathbb R_+$
the Green function
of the heat equation in $\mathbb R_+$ with Neumann boundary conditions at $0$, namely
   \begin{equation}
\label{3.5}
G_t^{\rm neum}(r,r') = G_t(r,r') +
G_t(r,-r'),\quad   G_t(r,r') =  \frac{e^{- \frac{(r-r')^2}{2t}}}{\sqrt{2\pi t}}
   \end{equation}
Observing that $G_t^{\rm neum}(r,r') =G_t^{\rm neum}(r',r)$ we write for $u\in \mathcal U$:
\[
G_t^{\rm neum}*u (r)= \int_{\mathbb R_+} G_t^{\rm neum}(r,r') u(r')\,dr'=
\int_{\mathbb R_+} G_t^{\rm neum}(r',r) u(r')\,dr'
\]

\end{definition}

\medskip

\begin{definition} (Barriers).
\label{defin:barriers}
Let $u\in  L^\infty(\mathbb R_+,\mathbb R_+)
\cap L^1(\mathbb R_+,\mathbb R_+)$ 
and such that $\dis{F(0;u )>0}$.
Then for all $\delta$ small enough   $u\in \mathcal U_\delta$
and for such $\delta$ we define
the ``barriers''  ${S_{k\delta}^{(\delta,\pm)}(u)}
$,  $k\in \mathbb N$,
as follows: we set
$S^{(\delta,\pm)}_{0}(u)=u$, and, for $k\ge 1$,
\begin{eqnarray}
\label{3.6}
&&
S^{(\delta,-)}_{k\delta}(u)= K^{(\delta)} G_\delta^{\rm neum} *S^{(\delta,-)}_{(k-1)\delta}(u)
\\&&
S^{(\delta,+)}_{k\delta}(u)= G_\delta^{\rm neum} * K^{(\delta)} S^{(\delta,{+})}_{(k-1)\delta}(u)
\nn
\end{eqnarray}

\end{definition}

\medskip
The functions $S_{k\delta}^{(\delta,\pm)}(u)$ deserve the name of ``lower and upper
barriers'' because we shall see that they form separated classes  with respect to the following notion of
partial order:

   \begin{definition}
\label{defin:ee8.1}
 (Partial order).
For any $u,v \in \mathcal U$
we set
   \begin{equation}
\label{4.1a}
u\le v \;\;\;\text{iff}\;\;\;
F(r; u)    \le F(r; v) \;\;\;\text{for all $r\ge 0$}
   \end{equation}
with  $F(r;u)$ as in \eqref{3.2}.
\end{definition}

When $u$ and $v$ have the same total mass then $u\le v$ if and only if $v$ can be obtained from $u$ by moving mass to  the right, this statement will be made precise in Proposition \ref{prop4.1}.

\medskip
{
In Proposition \ref{propineq}  we will prove many properties of the barriers and  in particular that they are monotone (namely the lower barrier is non increasing and the upper one is non decreasing) and that for all $\delta>0$
	$$
S_{k\delta}^{(\delta,-)}(u)\le S_{k\delta}^{(\delta,+)}(u)
	$$
Thus upper and lower barriers form separated classes and in the following Theorem we prove that there is a unique
separating element.

\medskip
\noindent
     \begin{theorem}[Barriers and separating elements]
           \label{thmineq}
Let $u \in L^\infty(\mathbb R_+,\mathbb R_+)\cap L^1(\mathbb R_+,\mathbb R_+)$ with $R(u) <\infty$, then
there exists a unique function $S_t(u)(r)$ continuous in $(r,t)$ for $t>0$
such that for all $t>0$ and $r\in \mathbb R_+$:
    \begin{eqnarray}
&&
{F(r;S_{t}(u)) = \lim_{\ell\to \infty} F(r;S_{t}^{(2^{-\ell}t,\pm)}(u))}\quad {\rm monotonically}
\label{4.4a}\\&&
F(r;S_{t}(u)) =\inf_{\delta: t=k\delta, k\in \mathbb N} F(r;S_{t}^{(\delta,+)}(u))
=\sup_{\delta: t=k\delta, k\in \mathbb N} F(r;S_{t}^{(\delta,-)}(u))
\label{4.4b}
  \end{eqnarray}
Moreover $S_{t}(u) \to u$ weakly as $t\to 0$
and if $u\le v$ 
then
$S_t(u) \le S_t(v)$.

\end{theorem}

\medskip

We have more detailed properties on the structure of
the barriers which are stated in Section \ref{sec:3}.

A key result of this paper is the identification
of the separating element between barriers
with the generalized solution of the FBP \ref{defin1.1}, {namely with the function $\bar \rho$ of Theorem \ref{exist}}.

\begin{theorem} [Characterization of ${\bar \rho}$]
\label{ineq}
The generalized solution $\bar \rho(\cdot,t)$
of the FBP \ref{defin1.1} with initial datum $\rho_0$  as in Definition \ref{defin1.2}
is equal to the
separating element $S_t(\rho_0)$ of Theorem \ref{thmineq}.
Moreover there is $a_2<\infty$ (which depends on $\rho_0$) so that
\cola{
$$
R(\bar\rho(\cdot,t)) <a_2 \qquad \text{for all} \quad t\ge 0,
$$
}
and
if $\rho_0$ is strictly positive in a neighborhood of the origin
then there exists   $a_1>0$
so that
\cola{
$$
R(\bar\rho(\cdot,t)) > a_1 \qquad \text{for all} \quad t\ge 0.
$$
}

\end{theorem}

\vskip.5cm

\subsection{\col{Scheme of the proofs}}
\col{The proofs of the results exploit analytical  and probabilistic arguments
and are organized in the following way.}

In Section \ref{sec:2} we prove   mass transport inequalities (some well known in the literature)
which are crucial in the subsequent steps.
These inequalities hold true with respect to the notion of partial order
in Definition \ref{defin:ee8.1} and with their help we deduce
in Section \ref{sec:3} properties of the  barriers which lead to the proof of
Theorem \ref{thmineq}.  In particular we
show that the sequence of upper and lower barriers
admits a unique limit that separates them.

\col{In Section \ref{sec:4} we recall the representation of the solution
of the heat equation in terms of Brownian motions and use it
to prove  part (a) of Theorem \ref{exist}, namely the
existence of an optimal sequence }{
of quasi-solutions.

\cola{In Section \ref{sec:5} we use the
uniqueness of the separating element between barriers}
to reduce the proofs of \cola{part (b), (c) of} Theorem
\ref{exist} and  \cola{Theorem \ref{ineq}}
to the variational problem of showing that
quasi-solutions are in between lower and upper barriers (modulo a ``small error'').
The crucial step is
Proposition \ref{keyes} in
Section \ref{sec:5.1}  where we prove that any
quasi-solution $(X_t,\rho(\cdot,t), {\eps})$, $t\in [0,T]$ is (up to a ``small error'') in between the barriers that start at $t=0$ from $\rho(\cdot,0)$.
As a 
\cola{consequence} of this result we prove in Section \ref{sec:8a} that:
\begin{itemize}
\item  any optimal sequence \cola{$(X^{(n)}_t,\rho^{(n)}(\cdot,t),\eps_n)$, $t\in (0,T)$} of quasi-solutions with initial datum $\rho_0$
is such that
the sequence of measures $\rho^{(n)}(r,t) dr$ on $\mathbb R_+$ is
tight;
\item any weak limit \cola{of $\rho^{(n)}(\cdot,t)$} is in between the barriers $S_{t}^{(2^{-\ell}t,\pm)}(\rho_0)$;
\item  by letting $\ell \to \infty$, and using Theorem \ref{thmineq} we then conclude that any weak limit is equal to $S_{t}(\rho_0)$, thus getting the statements (b) and (c) of Theorem \ref{exist};
\item the proof of Theorem \ref{ineq} is then a corollary of all this and it is given at the end of Section \ref{sec:8a}.
\end{itemize}

}

%
%
%

\vskip.5cm

\subsection{Remarks}

\vskip.5cm

\noindent
{\it {Generalized solutions}}.  We thus have global
existence of generalized solutions, hence a way
to continue
classical solutions past their singularity times (if they exist).
However we can only say that generalized solutions $\rho(r,t)$
are $(r,t)$ continuous and we do not know much about the motion of the edge.
In the existence part of the proof of Theorem \ref{exist}
we  construct a sequence of quasi-solutions with  the
edge which moves piecewise linearly (further work would be required to smoothen out
the discontinuities of its velocity) but we do not know about the regularity of the  motion
of the edge in the limit.
In Theorem \ref{ineq} we
prove
that at all times
\col{the edge} stays strictly positive and  it does not drift away to infinity.
However, if the edge and $\rho(r,t)$ are
``regular'' in some time interval then by Theorem \ref{ineq}
in that interval they define a classical solution.

\medskip
\vskip.5cm
\noindent
{\it {FBP with mass conservation}}.
The physics behind our problem is not the same as in the Stefan problem
where the existence of a boundary is related to an interface between two
phases.  Here instead it comes from the requirement of mass balance.
A more general formulation of the FBP with mass conservation would be the following.
Find pairs $(X_t,\rho(r,t))$,  where $t\ge 0$, $X_t>0$ and  $\rho(r,t)$ is for
each $t\ge 0$  a non negative function
which solves the  problem
\begin{eqnarray}
\label{12345}
&&\frac{\partial \rho}{\partial t}(r,t)   = \frac 12 \frac{\partial^2 \rho}{\partial r^2}(r,t)
 +f(r),\qquad  r \in(0,X_t)\nn\\
 &&\\&&
\rho( X_t,t)=0,\qquad \frac {\partial \rho}{\partial r}(0,t)=0,\qquad  \int_{0}^{X_t}\rho(r,t)\,dr =  \int_{0}^{X_0}\rho(r,0)\,dr\nn
\end{eqnarray}
where $f(r)$ is a non negative function with support in $[0,R_0]$,
$R_0<X_0$. The term $f(r)$ describes a mass source
of intensity $\int f(r) dr=j$,  with $j$ a positive constant, the mass then
diffuses according to the linear heat equation being reflected at
$0$ due to the Neumann condition, so that it can only escape
from $X_t$ where we impose Dirichlet condition. The problem is to determine
$X_t$ in such a way that the total mass in conserved.

\medskip

\noindent
The FBP \ref{defin1.1} fits within this scheme by choosing
$f=j D_0$, where $D_0$ denotes the Dirac delta at 0.  However our analysis
does not extend, at least directly, to smooth $f$, as in the mass transport inequalities
we exploit the fact that mass is added at the leftmost point.

\medskip
\vskip.5cm
\noindent
{\it {Scaling limits}}.
The FBP  \ref{defin1.1} (or similar versions) appeared in the study of the asymptotic
behavior, i.e. the scaling limit, of many different models in the context of statistical mechanics.
A non-exhaustive list of examples includes:

\vskip.3cm
\noindent
{\it {Particle systems}}.
The FBP  \ref{defin1.1}  has  been studied in connection
with the Fourier law in a moving domain where $\rho$ is interpreted as an
energy density, the sources are reservoirs which add and subtract the same
amount of energy at the boundaries while in the bulk energy diffuses
according to the linear heat equation. This is discussed in \cite{CDGP}
where the above FBP appears
as a heuristic guess for the continuum limit of a system of particles.
The particles move as independent, symmetric random walks on $\mathbb N$ with reflections at 0,
new particles are created at rate $j$
at the origin  while the rightmost particle is killed also at rate $j$.
To see the reason why the Dirichlet condition
appearing in \eqref{2.4c}
corresponds to the killing of the rightmost particle
we define $u(r,t):= \rho(r,t)$ for $r\in [0,X_t]$
and $u(r,t)=0$ elsewhere.  Then, in the weak formulation,
\begin{equation}
\label{1.4nnnn}
\frac{\partial u}{\partial t}   = \frac 12 \frac{\partial^2 u}{\partial r^2}
+ j D_0
- j D_{X_t}
,\qquad r\in \mathbb R_{+}
\end{equation}
where besides the source at the origin there is also a negative source at $X_t$
(the macroscopic counterpart of the microscopic killing at rate $j$).

Other particles systems whose scaling limit is given by a FBP of the Stefan type
have been studied in \cite{CS} and in \cite{dfp}, while -- in the reverse direction --
reference \cite{LV} studied a microscopic particle model for the Stefan freezing/melting
problem.

\medskip
\vskip.3cm
\noindent
{\it {Polymers}}.
In \cite{lacoin} a two-sided version (i.e. with a left and a right free boundaries)
of  the FBP \ref{defin1.1} in the interval $[-1,1]$ without the
condition \eqref{2.4b} (\col{that yields mass conservation}) has been studied.
The FBP is proved to be the scaling limit (under diffusive scaling) of a random dynamics
of polymers with pinning to a substrate.
However this problem is substantially different from ours, since mass
is not conserved anymore and it is decreasing with time.
As a consequence, for Lipschitz initial datum, existence and regularity
of the classical solutions could be proved until a maximal time where the two boundaries
collide (the proof uses the approach developed in \cite{CK}).

\medskip
\vskip.3cm
\noindent
{\it {Queuing theory}}.
The FBP  \ref{defin1.1}  has also natural connections
to fluid-limits in queuing theory \cite{SRI,ANS}. It is well known (Burke theorem)
that a tandem Jackson network with customers arriving
at the first queue as a Poisson process of intensity $\lambda >0$
gives rise to an outgoing flow of served clients at the
last queue that is also Poisson distributed with the same intensity.
The hydrodynamic limit of such tandem Jackson network
is given by the linear heat equation.
The FBP problem studied in this paper is a candidate for the scaling
limit of a tandem Jackson network of variable length.


\vskip1cm

\section{Mass transport inequalities}
\label{sec:2}

For the reader's convenience we present in this section some well known
facts about mass transport and use them to prove some properties which will then be
extensively used in the sequel. We start with  some elementary properties of the cut-and-paste operator
$K^{(\delta)}$ and the diffusion kernel $G_\delta^{\rm neum}$.

\medskip

\begin{proposition}[Preliminaries]
\label{prop3.1}

Let $\delta>0$ and let $u$ and $v$ be in $\mathcal U_\delta$, then
   \begin{equation}
\label{3.7}
F(0; u) = F(0;K^{(\delta)} u)=F(0;G^{\rm neum}_\delta *u)
   \end{equation}
      \begin{equation}
\label{3.8}
|K^{(\delta)} u -K^{(\delta)} v|_1\le |u-v|_1,\quad
|G^{\rm neum}_\delta *u-G^{\rm neum}_\delta *v|_1\le |u-v|_1
   \end{equation}
   \begin{equation}
	\label{3.9}
|(K^{(\delta)}-1)u|_1 {=} 2j\delta
\end{equation}
   \begin{equation}
	\label{3.10}
|S^{(\delta,\pm)}_{k\delta}(u) - S^{(\delta,\pm)}_{k\delta}(v)|_1 \le |u-v|_1,\qquad \forall k\in \mathbb N
\end{equation}
%
%

\end{proposition}

\medskip

\noindent
{\bf Proof.} \eqref{3.7} follows directly from the definition of
$K^{(\delta)}$ and $G_\delta^{\rm neum}$.  To prove the first inequality
in \eqref{3.8}  we write $u=c_uD_0+ \rho_u$, $v= c_vD_0+\rho_v$ and
assuming that $R_\delta(u)\le R_\delta(v)$,
      \begin{equation*}
|K^{(\delta)}u -K^{(\delta)}v|_1= |c_u-c_v| + \int_0^{R_\delta(u)} |\rho_u-\rho_v|
+ \int_{R_\delta(u)}^{R_\delta(v)} \rho_v
   \end{equation*}
We can then add $\dis{\int_{R_\delta(v)}^\infty \rho_v}$ and subtract
   $\dis{\int_{R_\delta(u)}^\infty \rho_u}$ as they are both equal to $j\delta$:
      \begin{equation*}
|K^{(\delta)} u -K^{(\delta)} v|_1= |c_u-c_v| + \int_0^{R_\delta(u)} |\rho_u-\rho_v|
+ \int_{R_\delta(u)}^{\infty} (\rho_v-\rho_u) \le |u-v|_1
   \end{equation*}
Analogously
      \begin{eqnarray*}
|G^{\rm neum}_\delta *u-G^{\rm neum}_\delta *v|_1 \le |c_u-c_v| \int  G^{^{\rm neum}}_\delta(0,r')\,dr'
+\int\int |\rho_u(r)-\rho_v(r)|G^{^{\rm neum}}_\delta(r,r')\,dr\,dr' &&\\
{=|c_u-c_v| +\int |\rho_u(r)-\rho_v(r)|\,dr} &&
   \end{eqnarray*}
which is equal to $|u-v|_1$.  Finally
   \begin{equation*}
|(K^{(\delta)}-1)u|_1 {=} j\delta + \int_{R_{\delta}(u)}^\infty \rho_u = 2j\delta
\end{equation*}
while \eqref{3.10} is a direct consequence of \eqref{3.8}.
\qed


\vskip.5cm

{\begin{proposition}[Mass displacement]
\label{prop4.1}
Given  $u\le v$ in $ L^\infty \cap L^1$  with $F(0;u)=F(0;v)$ we define for $r\in \mathbb R_+$:
  \begin{equation}
\label{4.5a}
f(r) := \sup\;\Big\{r': \int_0^{r'} v(z)\,dz= \int_0^r u(z)\,dz\Big\},
   \end{equation}
Then
  \begin{equation}
\label{4.6a}
f(r) \ge r
   \end{equation}
and for any function $\phi\in L^\infty(\mathbb R_+,\mathbb R)$
  \begin{equation}
\label{4.7a}
 \int_0^\infty v(r)\phi(r)\,dr = \int_0^\infty u(r) \phi(f(r))\,dr
   \end{equation}
\end{proposition}}

\medskip

\noindent
{\bf Proof.}  Since $F(0;u)= F(0;v)$,
\[
\int_0^r u(z) \,dz + F(r;u)= \int_0^r v(z) \,dz + F(r;v)
\]
and since $F(r;u)\le F(r;v)$
\[
\int_0^r u(z) \,dz  \ge \int_0^r v(z) \,dz
\]
which yields \eqref{4.6a}.
By a density argument \eqref{4.7a} follows from \eqref{4.6a}.


\qed

%
%
%
\medskip

\begin{corollary}[Moving mass to the right]
Let  $u\le v$ in $ L^\infty \cap L^1$ and $F(0;u)=F(0;v)$, then for all bounded, non decreasing
functions $h$ on $\mathbb R_+$:
  \begin{equation}
\label{4.7aaa}
  \int_0^\infty u(r) h(r) \,dr  \le \int_0^\infty v(r) h(r)\,dr
   \end{equation}

\end{corollary}

\medskip

\noindent
{\bf Proof.}  Observe that \eqref{4.7aaa} is verified by definition for all
functions $h$ of the form $\mathbf 1_{[R,\infty)}$, $R\ge 0$.  Its validity for
functions $h$ as in the text follows from \eqref{4.7a} because
  \begin{equation*}
 \int_0^\infty v(r)h(r)\,dr =  \int_0^\infty u(r) h(f(r)) \,dr
   \end{equation*}
and $ h(f(r)) \ge h(r)$ by \eqref{4.6a}.  \qed

\vskip.5cm

\begin{lemma}[Left cut]
\label{lemma4.2}
Let $u\le v$ and assume that $m:= F(0;v)-F(0;u)>0$
Define $\tilde R$  so that $\dis{\int_0^{\tilde R} v= m}$,
then
  \begin{equation}
\label{4.8a}
u\le  v \;{\mathbf 1_{[\tilde R,+\infty)}}=:\tilde v,\quad
F(0;u)=F(0;\tilde v)
   \end{equation}


\end{lemma}
\medskip

\noindent
{\bf Proof.}   From  $F(0;u) = F(0;v)-m=F(0,\tilde v)$ we get
	\begin{eqnarray*}
	\int_0^r u(z)dz + F(r;u) = \int_0^r v(z)dz + F(r;v) -m
	\end{eqnarray*}
Since  $F(r;u)\le F(r;v)$, for all  $r\ge \tilde R$
	\begin{eqnarray*}
	\int_0^r u(z)dz\ge \int_0^r v(z)dz - m=\int_{\tilde R}^r v(z)dz=\int_0^r \tilde v(z)dz
	\end{eqnarray*}
Also for $r<\tilde R$ we have
$\dis{\int_0^r u\ge 0= \int_0^r \tilde v}$, so that $\dis{\int_0^r u\ge  \int_0^r \tilde v}$
for all $r$.
Since $F(0;u) = F(0;\tilde v)$ the previous inequality implies that
$F(r;u) \le F(r;\tilde v)$.
\qed

\medskip

\noindent
\begin{definition}[Partial order modulo $m$]
\label{def:2.1}
For any $u$ and $v$ in $\mathcal U$ and  $m>0$, we define  \begin{equation}
\label{4.9a}
u\le v\;\;\text{ modulo $m$ \; iff \;\;for all $r\ge 0$:}\;\;  F(r;u) \le F(r;v) + m
   \end{equation}

\end{definition}

\vskip.5cm

\begin{lemma}[Right cut]
\label{prop4.2}

Let $u\le v$ modulo $m$, then

\medskip
  \begin{equation}
\label{4.10a}
u^* := u \;{\mathbf 1_{[0,R_m]}}\le v \;\;\;\;\text{with}\;\; R_m = \inf\{ r: F(r;u) = m\}
   \end{equation}

\end{lemma}

\medskip

\noindent
{\bf Proof.} For $r\le R_m$, $F(r;u^*)= F(r;u)-m$ so that $F(r;u^*) \le F(r;v)$.
For $r\ge R_m$, $F(r;u^*)=0 \le F(r;v)$.  \qed


\medskip

 \begin{lemma}[Partial order under diffusion and cut-and-paste]
 \label{lemma4.5}
Let $u\le v$ in $\mathcal U$, then for any $t>0$
  \begin{equation}
\label{4.11a}
G_t^{\rm neum} *u \le G_t^{\rm neum} *v
   \end{equation}
and if $u\le v$ with $u,v$ in $\mathcal U_\delta$ then
  \begin{equation}
\label{4.12a}
K^{(\delta)} u \le u,\quad
K^{(\delta)} u \le K^{(\delta)} v,\quad u \; {\mathbf 1_{[0,R_\delta(u)]}} \le v
{\;\mathbf 1_{[0,R_\delta(v)]}}
   \end{equation}
 \end{lemma}

\medskip

\noindent
{\bf Proof.}  It is clearly enough to prove {$G_t^{\rm neum} *u \le G_t^{\rm neum} * \tilde v
$ with $\tilde v$ as in \eqref{4.8a}.  We have   $u \le \tilde v$, $F(0;u)=F(0,\tilde v)$, then,  if} $f$ is the map defined in
Proposition \ref{prop4.1} relative to $u$ and $\tilde v$, by \eqref{4.7a}
 with
\[
\phi(r)= \int_{R}^\infty G_t^{\rm neum}(r,r')\,dr'
\]
we get
  \begin{equation*}
F(R;G_t^{\rm neum}* \tilde v)-F(R;G_t^{\rm neum}* u)=\int u(r) \Big[\phi(f(r))- \phi(r)\Big] \,dr
   \end{equation*}
By an explicit computation: $\dis{ \frac{d}{dr}\int_{R}^\infty G_t^{\rm neum}(r,r')\,dr' >0}$, moreover,  by Proposition \ref{prop4.1}, $f(r) \ge r$, then
 \eqref{4.11a} follows {because $\phi(f(r))\ge \phi(r)$}.

The inequality $K^{(\delta)}u \le u$ trivially follows from the definition \eqref{3.4}.
Furthermore  we have
      \begin{equation*}
K^{(\delta)}u -K^{(\delta)} v= (c_u-c_v)D_0 +  (\rho_u-\rho_v) \;{\mathbf 1_{[0,R_\delta(u)]}}
-\rho_v \, {\mathbf 1_{(R_\delta(u), R_\delta(v)]}}
   \end{equation*}
where $R_\delta(\cdot)$ is defined in \eqref{3.4}, {and $R_\delta(u)\le R_\delta(v)$ because $u\le v$.}  Hence
      \begin{equation*}
F(r;K^{(\delta)} u)-F(r; K^{(\delta)} v)= \Big(F(r;u)-F(r;v)\Big){\mathbf 1_{[0,R_\delta(u)]}} -{\mathbf 1_{(R_\delta(u), R_\delta(v)]}} \int_r^{R_\delta(v)}
\rho_v(r') dr'
   \end{equation*}
which is therefore $\le 0$. Same argument proves the last inequality in \eqref{4.12a}.

 \qed

\medskip

 \begin{lemma}[Partial order modulo $m$ under diffusion and cut-and-paste]
Let $u\le v$ modulo $m$, $u,v$ in $\mathcal U_\delta$, then
  \begin{equation}
\label{4.13a}
 G^{\rm neum}_t *u \le  G^{\rm neum}_t *v\;\;\; \text{modulo }\;\; m
   \end{equation}
  \begin{equation}
\label{4.14a}
 u \le K^{(\delta)}   v\;\;\; \text{modulo }\;\; m+j\delta
   \end{equation}
and if $F(0;u) \le F(0;v)+\alpha $ and $m\ge j\delta$, $\alpha\ge 0$, then
  \begin{equation}
\label{4.15a}
 K^{(\delta)}  u \le  v\;\;\; \text{modulo}\;\;\; \max (m -j\delta, \alpha)
   \end{equation}

 \end{lemma}

\medskip

\noindent
{\bf Proof.}

{\em Proof of \eqref{4.13a}.}
By \eqref{4.10a} $u^* \le v$ and $F(0;u-u^*)=m$.  Then by \eqref{4.11a}
\begin{eqnarray*}
F(r;G^{\rm neum}_t*u) &=& F(r;G^{\rm neum}_t *u^*) + F(r;G^{\rm neum}_t *(u-u^*))
\\&\le& F(r;G^{\rm neum}_t *v) + F(0;G^{\rm neum}_t *(u-u^*)) = F(r;G^{\rm neum}_t *v) +m
\end{eqnarray*}

\medskip
 {\em Proof of \eqref{4.14a}.}  We need to prove that for all $r\ge 0$:
\[
F(r;u) \le F(r;K^{(\delta)} v) + m +j\delta
\]
For $r=0$ we have $F(0;v)=F(0; K^{(\delta)}v)$ so that
\[
F(0;u) \le F(0;v)+m =F(0; K^{(\delta)}v)+m
\]
We next consider
$0<r\le R_\delta(v)$ for which we have
  \begin{equation*}
 F(r; K^{(\delta)}   v) = F(r;v) - j\delta \ge F(r;u) - m -j\delta
   \end{equation*}
while for  $r> R_\delta(v)$ then {$F(r;v)\le j\delta$ hence}
  \begin{equation*}
 F(r; K^{(\delta)}   v) = 0 \ge  F(r;v) - j\delta \ge F(r;u) - m -j\delta
   \end{equation*}
hence \eqref{4.14a}.

{\em Proof of  \eqref{4.15a}. }  We first observe that for $r=0$
\[
F(0; K^{(\delta)} u) = F(0;u) \le F(0;v)  +\alpha
\]
As before for $0<r< R_\delta(u)$,
\[
F(r; K^{(\delta)} u)=  F(r;u)  -j\delta \le F(r;v)+m-j\delta
\]
For $r\ge R_\delta(u)$, $F(r;K^{(\delta)} u)=0 \le F(r;v) \le  F(r;v) + m-j\delta$.
\qed

\vskip.5cm


\section{Proof of Theorem \ref{thmineq}}
\label{sec:3}
We first state and prove properties of the barriers which will be often used in the sequel.
These properties
have been proved in \cite{CDGP} when $G_t^{\rm neum}(r,r')$ is
replaced by the Green function with Neumann condition both at 0 and at 1. The proofs
in \cite{CDGP} extend to our case but the extension is not  straightforward because
in \cite{CDGP} properties of the particles evolution are sometimes used.
We shall give in this paper a self contained proof which is purely
analytical and uses only elementary properties of mass transport inequalities.

\medskip

\begin{proposition} [Barrier properties]
\label {propineq}

Let  {$u,v \in L^\infty(\mathbb R_+,\mathbb R_+)\cap L^1(\mathbb R_+,\mathbb R_+)$} 
and
such that $j\delta<F(0;u)$, {$j\delta<F(0;v)$}.  Then

\begin{itemize}

\item If $u\le v$ then for all $k\in \N$
\begin{equation}
\label{uu1}S_{k\delta}^{(\delta,\pm)}(u)\le S_{k\delta}^{(\delta,\pm)}(v)\end{equation}

\item  Let $\delta = \ell\delta'$ with  $\ell$ a positive integer, then for all $k\in \N$

  \begin{eqnarray}
\label{4.3a.0}
&&S_{k\delta}^{(\delta,-)}(u) \le S_{k\delta}^{(\delta',-)}(u) \le S_{k\delta}^{(\delta',+)}(u)\le
S_{k\delta}^{(\delta,+)}(u)
\end{eqnarray}

\item  For all $k\in \N$
\begin{eqnarray}
&&
|S_{k\delta}^{(\delta,+)}(u) - S_{k\delta}^{(\delta,-)}(u) |_1\le 4j\delta
\label{4.3a.00}
  \end{eqnarray}

%
%

\end{itemize}

\end{proposition}

\medskip
\noindent
{\bf Proof.}

\begin{itemize}
\item
{\em Proof of \eqref{uu1}.} It follows immediately from Lemma \ref{lemma4.5}.

\vskip.3cm
\item {\em Proof of  \eqref{4.3a.0}.}  It is a consequence of
the following inequalities which will be  proved next, one after the other:
		\begin{equation}
		\label{nna.1}
S_{k\delta}^{(\delta,-)}(u) \le
S_{k\delta}^{(\delta,+)}(u)\qquad \forall k \in\mathbb N
			\end{equation}
		\begin{equation}
		\label{nna.2}
S_{k\delta}^{(\delta',+)}(u) \le
S_{k\delta}^{(\delta,+)}(u),\qquad S_{k\delta}^{(\delta,-)}(u) \le S_{k\delta}^{(\delta',-)}(u)\qquad \text{for $\delta = \ell \delta'$} \quad \forall \ell, k \in \mathbb N
			\end{equation}

\vskip.3cm
\begin{itemize}
\item[$\star$] {\em Proof of \eqref{nna.1}.} We first prove that 
\begin{equation}
		\label{nna}
S_{\delta}^{(\delta,-)}(u) \le
S_{\delta}^{(\delta,+)}(u)
			\end{equation}
By \eqref{4.14a}, $u \le K^{(\delta)} u$ modulo $j\delta$, so that using  \eqref{4.13a}
\[
G^{\rm neum}_\delta *u \le G^{\rm neum}_\delta *K^{(\delta)}u=S_{\delta}^{(\delta,+)}(u)\quad \text{modulo}\;\; j\delta
\]
Then, {since $F(0; G^{\rm neum}_\delta *u )=F(0; S_{\delta}^{(\delta,+)}(u))$,} by \eqref{4.15a} with $\alpha=0$ and $m=j\delta$ we get
$$
 S_{\delta}^{(\delta,-)}(u)=K^{(\delta)} G^{\rm neum}_\delta *u \le S_{\delta}^{(\delta,+)}(u)
$$
The inequality \eqref{nna} is then proved.
\vskip.2cm

%
%
%
We shall now prove  \eqref{nna.1} by induction on $k$.  The inequality for $k=1$
follows from \eqref{nna}. We next suppose
that it holds for $k$ and need to prove it for $k+1$.  Call
$u'= S_{k\delta}^{(\delta,-)}(u)$ and $v'=S_{k\delta}^{(\delta,+)}(u)$ then by
the induction assumption $u'\le v'$, so that, {from \eqref{uu1} and \eqref{nna},}
$$S_{\delta}^{(\delta,-)}(u') {\le S_{\delta}^{(\delta,-)}(v')} \le S_{\delta}^{(\delta,+)}(v')$$
hence  \eqref{nna.1} follows because
\[
S_{(k+1)\delta}^{(\delta,-)}(u)=S_{\delta}^{(\delta,-)}(u')\qquad
\text{and} \qquad S_{(k+1)\delta}^{(\delta,+)}(u)=S_{\delta}^{(\delta,+)}(v')
\]
\vskip.3cm
\item[$\star$] {\em Proof of \eqref{nna.2}.} We first prove it for $k=1$:
 \begin{equation}
\label{4.16a}
S_{\delta}^{(\delta',+)}(u)\le
S_{\delta}^{(\delta,+)}(u)\qquad \text{for  $\delta = h\delta'$, }\quad h \in \mathbb N
   \end{equation}
We have
{
\[
S_{\delta}^{(\delta',+)}(u) =S_{h\delta'}^{(\delta',+)}(u) = G^{\rm neum}_{\delta'} * K^{(\delta')} \cdots G^{\rm neum}_{\delta'} * K^{(\delta')} u \qquad \text{$h$ times}
\]}
\[
S_{\delta}^{(\delta,+)}(u) =S_{h\delta'}^{(\delta,+)}(u) =  G^{\rm neum}_{\delta'} *  \cdots *G^{\rm neum}_{\delta'} *  K^{(\delta)}u \qquad \,\,\,\text{$h$ times}
\]
{By applying \eqref{4.14a} and \eqref{3.7}
 we get
	\begin{eqnarray*}
&&u\le K^{(\delta)}(u)\quad  \text{modulo } j\delta \qquad \text{and} \quad  F(0; u)=F(0; K^{(\delta)}u)
\end{eqnarray*}
then, using \eqref{4.15a} with $\alpha =0$ and $m=j\delta$ we get
\begin{eqnarray*}
 K^{(\delta')} u\le
K^{(\delta)}u \quad  \text{modulo } \: j\delta-j\delta'
\end{eqnarray*}
then, from \eqref{4.13a},
\begin{eqnarray*}
G_{\delta'}^{\rm neum}*K^{(\delta')} u\le
G_{\delta'}^{\rm neum}*K^{(\delta)}u \quad  \text{modulo } \: j\delta-j\delta'
\end{eqnarray*}
Let now $u':=G_{\delta'}^{\rm neum}*K^{(\delta')} u$ and $v':=G_{\delta'}^{\rm neum}*K^{(\delta)}u$ then we have $u'\le v'$ modulo  $j\delta-j\delta'$ and $F(0;u')=F(0;v')$. Thus we can apply again \eqref{4.15a} with $\alpha =0$ and $m=j\delta-j\delta'$ getting $K^{(\delta')} u'\le v'$ modulo $j\delta-2 j\delta'$. From  \eqref{4.13a} we then get
\begin{eqnarray*}
G_{\delta'}^{\rm neum}*K^{(\delta')} u'\le
G_{\delta'}^{\rm neum}*v' \quad  \text{modulo } \: j\delta-2 j\delta'
\end{eqnarray*}
 Then we obtain
 \eqref{4.16a} by iteration.}

We prove now  \eqref{nna.2} by  induction on $k$. We have just proved it for $k=1$,
suppose then that
it is verified for $k$. Call $u'= S_{k\delta}^{(\delta',+)}(u)$, $v'=S_{k\delta}^{(\delta,+)}(u)$
and use \eqref{4.16a} and the induction assumption $u'\le v'$.  Then, {from \eqref{uu1} and \eqref{4.16a} we get}
\[ S_{(k+1)\delta}^{(\delta',+)}(u) = S_{\delta}^{(\delta',+)}(u')
\le {S_{\delta}^{(\delta',+)}(v') \le S_{\delta}^{(\delta,+)}(v') =  S_{(k+1)\delta}^{(\delta,+)}(u)}
\]
The proof that $S_{k\delta}^{(\delta,-)}(u) \le S_{k\delta}^{(\delta',-)}(u)$
is similar and omitted.

\end{itemize}
 \eqref{4.3a.0} is then proved.
\vskip.3cm

\item {\em Proof of  \eqref{4.3a.00}.}
Shorthand $G$ for the operator $G_\delta^{\rm neum} *$ and
 \[
\phi:= K^{(\delta)}G \cdots  K^{(\delta)}G u,\qquad \psi:=G  K^{(\delta)}\cdots G K^{(\delta)} u \qquad {\text{$k$ times}}
 \]
so that we need to bound the total variation of $\phi-\psi$. Call
 \[
 v = K^{(\delta)}u,\qquad v_k= G K^{(\delta)}\cdots G v,\qquad u_k=G K^{(\delta)}\cdots G u \qquad {\text{$k$ times}}
 \]
Thus $u_k$ and $v_k$ are obtained by applying $G (K^{(\delta)}G)^{k-1}$ to $u$ and respectively $v$, {hence $\phi= K^{(\delta)}u_k$ and $\psi=v_k$.}
{From \eqref{3.8} we have that $G (K^{(\delta)}G)^{k-1}$ is a contraction, then,} using  \eqref{3.9}, we  get
 \begin{eqnarray*}
 |\psi-\phi|_1 &{=}& |K^{(\delta)}u_k - v_k|_1 \le |K^{(\delta)}u_k - u_k|_1 +
|v_k-u_k|_1\\
  &{=}& 2j\delta + |v_k-u_k|_1 \le 2j\delta +|u-v|_1 {=} 4j\delta
 \end{eqnarray*}

\end{itemize}

\vskip.5cm
\noindent
{\bf   Proof of Theorem \ref{thmineq}.}
We start by proving the existence of $S_t(u)$.
Fix $\tau>0$ and for any $\ell\in \mathbb N$ define $u_{\ell;\tau}(r,t)$ as the
linear interpolation of the functions $S^{(2^{-\ell}\tau,+)}_{k 2^{-\ell}\tau}(u)$, $k\in \mathbb N$.
In \cite{CDGP} it is proved that restricted to $t\ge \si>0$ the family
$\{u_{\ell;\tau}\}$ is equibounded and equicontinuous.  The proof uses only properties
of the Green function which are valid both in the domain $[0,1]$ considered in \cite{CDGP}
and in our domain $[0,\infty)$.  Thus referring to {Theorem 2.3} in \cite{CDGP}
we can say (by the Ascoli-Arzel\`a theorem and a diagonalization procedure)
that  for $t>0$, $u_{\ell;\tau}(r,t)$ converges by subsequences as $\ell\to \infty$
to a continuous function $u_{\tau}(r,t)$ which in principle depends also on the
converging subsequence.  We now prove that for all $r\ge 0$
and $t\in \{ k 2^{-\ell}\tau; \;\;k,\ell \in \mathbb N_+\}$,
  \begin{equation}
\label{4.17a}
\lim_{\ell\to \infty} F\Big(r;S_{t}^{(2^{-\ell}\tau,+)}(u)\Big)=F\big(r;u_{\tau}(\cdot,t)\big)
   \end{equation}
By Fatou's lemma 
     \begin{equation*}
\lim_{\ell\to \infty} F\Big(r;S_{t}^{(2^{-\ell}\tau,+)}(u)\Big)\ge F\big(r;u_{\tau}(\cdot,t)\big)
   \end{equation*}
so that we just need to prove the converse inequality.
Let $R>r$, then
\[
F\Big(r;S_{t}^{(2^{-\ell}\tau,+)}(u)\Big) = \int_r^R S_{t}^{(2^{-\ell}\tau,+)}(u) + \int_R^\infty S_{t}^{(2^{-\ell}\tau,+)}(u)
\]
The first integral converges to $\dis{\int_r^R u_{\tau}(\cdot,t)}$ by the Lebesgue
dominated convergence theorem while by \eqref{4.3a.0}, \eqref{4.12a} {and \eqref{4.11a}}
\[
\int_R^\infty S_{t}^{(2^{-\ell}\tau,+)}(u) \le \int_R^\infty S_{t}^{(\tau,+)}(u) \le
\int_R^\infty G_t^{\rm neum}*u
\]
which decays exponentially as $R\to \infty$ (recall that $u$ has compact support).
This proves \eqref{4.17a}.

From \eqref{4.17a} it follows that the limit $u_\tau$
is independent of the subsequence,
hence that $u_{\ell;\tau}(r,t)$ converges  as $\ell\to \infty$ and for $t>0$
to  $u_{\tau}(r,t)$.  In \cite{CDGP}, see {Theorem 8.4}, it is proved that
$u_{\tau}(r,t)$ is actually independent of $\tau$,
the proof extends straightforwardly to our case and it
is omitted. We can then identify  $S_{t}(u) : = u_{\tau}(\cdot,t)$ and by its definition
and \eqref{4.3a.0} we know that {the equality for the upper barrier   in \eqref{4.4a} is satisfied, i.e.
$\dis{F(r;S_{t}(u)) =\lim_{\ell\to \infty} F(r;S_{t}^{(2^{-\ell}\tau,+)}(u) )}$. Using \eqref{4.3a.00} we get
\begin{equation}
\Big|F(r;S_{t}^{(2^{-\ell}\tau,-)}(u))-F(r;S_{t}(u))\Big| \le 4j 2^{-\ell}\tau + \Big|F(r;S_{t}^{(2^{-\ell}\tau,+)}(u))-F(r;S_{t}(u))\Big|
\end{equation}
then  the equality in \eqref{4.4a}  is satisfied also for the lower barrier.
 From \eqref{4.3a.0} and \eqref{4.4a} it necessarly follows that}
    \begin{equation}
\label{4.433a}
F(r;S_{t}^{(2^{-\ell}\tau,-)}(u) ) \le
 F(r;S_{t}(u)) \le F(r;S_{t}^{(2^{-\ell}\tau,+)}(u) )
  \end{equation}
this and \eqref{4.4a} yield \eqref{4.4b}.

The proof that $S_t(u) \to u$ weakly as $t\to 0$ is essentially the same as
in \cite{CDGP}, Proposition 8.1, and it is omitted.  The last statement in the theorem,
namely the fact that if $u\le v$
then  $S_t(u) \le S_t(v)$, follows {directly  from \eqref{4.4a} and \eqref{uu1}.}
%
%
\qed

\vskip.5cm


\section{\col{Existence of quasi-solutions with arbitrary good accuracy}}
\label{sec:4}

\medskip
The analysis in this section, as well as in the following one, exploits extensively
a probabilistic representation of the solutions of
the heat equation in terms of  Brownian motions.
We first recall  the basic representation formulas and refer to the literature
for their proofs.
At the end of the section  we prove part (a) of  Theorem \ref{exist}, namely
the existence of \cola{optimal} quasi-solutions.

%
\vskip.5cm

\subsection{Probabilistic representation of the Green functions}

Let $P_{r;s}$, $r\ge 0$, $s\ge 0$ be the law of the
Brownian motion $B_t$, $t\ge s$, which starts from $r$ at time $s$, i.e.\  $B_s=r$, and which is
reflected at $0$.
The law of $B_t$ is absolutely continuous with respect to the
Lebesgue measure and has a probability density  that we
denote by $G^{\rm neum}_{s,t}(r,r')$. As a result
\begin{equation}
\label{3.1}
\rho(r,t) := \int G^{{\rm neum}}_{{s},t}(r',r) \rho(r',s) \,dr',\;\;\;\qquad G^{{\rm neum}}_{s,t}(r',r)dr=P_{{r'};s}[B_t\in (r,r+dr)]
\end{equation}
is the solution of the heat equation in $\mathbb R_+$ with
Neumann conditions at 0 and  datum $\rho(r',s)$ at time $s$; therefore
$G^{\rm neum}_{s,t}(r,r')$ is its corresponding Green function.

%
%

\vskip.5cm

%
%
%

Call $X=( X_t$, $t\ge 0)$ and denote for
$s\ge 0$ 
 \begin{equation}
\label{5.2.1}
 \tau^X_s = \inf\{t\ge s:  B_t \ge X_t\},\;\;\;\text{and $=\infty$ if the set is empty}
\end{equation}
The inf is a minimum as $X_t$  is continuous (see Definition \ref{defin1.2}).
The law of $B_t$, restricted to trajectories so that  $\{ \tau^X_s >t\}$,
has a density with respect to Lebesgue that we denote
by  $G^{X,\,{\rm neum}}_{s,t}(r,r')$: for any interval $I\subset \mathbb R_+$
 \begin{eqnarray}
\label{3.1444441111111}
&&
\int_I G^{X,\,{\rm neum}}_{s,t}(r',r)dr=P_{{r'};s}[\tau^X_s >t\,;\,B_t\in I]
\end{eqnarray}
If $\rho(r',s)\in L^\infty([0,X_s),\mathbb R_+)$, then the solution
of \eqref{2.4a}--\eqref{2.4c} (with initial datum   $\rho(r',s)$ at time $s$)
is  given by
\begin{eqnarray}
\label{3.144444}
&& \rho(r,t) := \int  G^{X,\,{\rm neum}}_{{s},t}(r',r) \rho(r',s) \,dr'
+ \int_s^t jG^{X,\,{\rm neum}}_{s',t}(0,r) \,ds',
\end{eqnarray}
and therefore
$G^{{X,}\rm neum}_{s,t}(r,r')$ is its corresponding Green function (recall that $X_t$
is piecewise $C^1$).

%
%

\vskip.5cm

We also have a nice representation for the mass
$ \Delta^X_{[0,t]}(u)$ which is removed from the system in the time interval
$[0,t]$
(due to the Dirichlet boundary conditions when the initial datum is $u$)
in terms of the probability that the Brownian motion reaches
the edge $X_t$.

\medskip
\begin{lemma}[Mass loss]  Let $\rho(r,t)$ solve \eqref{2.4a}--\eqref{2.4c} with initial datum $u\in L^\infty([0,X_0),\mathbb R_+)$ at time 0.  Then
\begin{equation}
\label{5.3.1.1}
F(0,\rho(\cdot,t)) =  F(0,u) + jt - \Delta^X_{[0,t]} (u)
\end{equation}
where
 \begin{eqnarray}
\label{5.3.1}
&& \Delta^X_{[0,t]} (u) =   \int   u(r') P_{r';0}\Big[\tau_0^X \le t\Big] \,dr'
 + \int_0^t jP_{0;s}\Big[ \tau^X_s\le t\Big] \,ds
\end{eqnarray}
In particular mass conservation as in  \eqref{1.4d}
requires that $\Delta^X_{[0,t]}(u) = jt$.
\end{lemma}

\medskip

\noindent
{\bf Proof.}
By integration of \eqref{3.144444}, the total mass at any time $t$ is given by
\begin{equation}\label{5.3.1.2.2}
F(0,\rho(\cdot,t)) = \int dr' {u(r')} P_{r',0}[\tau_0^X > t]
+ \int_{0}^{{t}} j ds P_{0,s}[\tau_s^X > t]\;
\end{equation}
Writing $P_{r',0}[\tau_0^X > t] = 1 - P_{r',0}[\tau_0^X \le t]$
and $P_{0,s}[\tau_s^X > t] = 1- P_{0,s}[\tau_s^X \le t]$
one finds
\begin{equation}\label{5.3.1.2}
F(0,\rho(\cdot,t)) =  F(0,{u}) - \int {u(r')} P_{r',0}[\tau_0^X \le t]dr'
+ j t - \int_{0}^{{t}} j  P_{0,s}[\tau_s^X \le t]ds\;.
\end{equation}
Then \eqref{5.3.1.1} follows from \eqref{5.3.1.2}.
\qed

\medskip

The  distribution of $\tau^X_0$
(inherited from the
probability measure
$P_{r;0}$) conditioned to the event  $\tau^X_0 \le t$ is denoted by
$\la^X_{r,t}(ds)$.
Hence  if $I$ is any interval in $\mathbb R_+$, {$t'\ge t$} then
   \begin{equation}
P_{r;0}\Big[ B_{{t'}} \in I \:\big|\, \tau^X_0 \le t \Big] = \int_0^t  P_{X_s,s}\big[ B_{{t'}} \in I]
\la^X_{r,t}(ds)
\label{7.18a}
	\end{equation}
Analogously denoting by  $\kappa^X_{s',t}$ the conditional probability density of $\tau^X_{s'}$
given that $\tau^X_{s'}\le t$ (under the probability measure
$P_{0;s'}$), we have that that {for any $t' \ge t$}
   \begin{equation}
    \label{7.19a}
 P_{0;s'}\Big[ B_{ {t'}} \in I \:\big|\, \tau^X_{s'} \le t \Big] = \int_{s'}^t P_{X_s,s}\big[ B_{{t'}} \in I]
 \kappa^X_{s',t}(ds).
	\end{equation}


\vskip.5cm

%
%

\subsection{Construction of quasi-solutions}

The following is a key lemma for the construction of quasi-solutions.

\medskip

\begin{lemma}
\label{ciccio}
Let $u\in L^\infty({[0,X_0]},\mathbb R_+)$,
$t^*>0$ and $X^{ V}_t=X_0+Vt$ with
$V> V^*:=- X_0/t^*$.  Call $u^{(V)}(r,t)$  the solution of \eqref{2.4a}--\eqref{2.4c} with initial datum   $u$ at time $0$
and edge $X_t^ V$.
Then there is a value of $V$ such that
\begin{equation}
\label{5.8888}
\Delta^{X^V}_{[0,t^*]}(u) = jt^*,\qquad \sup_{t\le t^*} |\Delta^{X^ V}_{[0,t]}(u) -jt| \le jt^*
\end{equation}

\end{lemma}

\medskip

\noindent
{\bf Proof.}
 The equality (on the left equation of
\eqref{5.8888}) follows directly from the following statements:
\begin{itemize}

\item[i)]  $\Delta^{X^ V}_{[0,t^*]}(u)$ converges to $jt^*+F(0;u)$ as $V\to V^*$
and it converges to 0 as $V\to \infty$.

\item[ii)] $\Delta^{X^ V}_{[0,t^*]}(u)$ depends continuously on $V$ in $(V^*,+\infty)$.
%
\end{itemize}

\noindent
The proof of such statements
is quite elementary as it can be reduced
to simple properties of the Green functions, for completeness we give some details.

\vskip.3cm
\begin{itemize}

\item[i)] 
Let $V>V^*$ and let $\eps:= X_0+ Vt^* = (V-V^*)t^*$.
Then by \eqref{5.3.1.1} and \eqref{5.3.1.2.2}
 \begin{eqnarray*}
&&\ 0\le jt^*+F(0;u) - \Delta^{X^ V}_{[0,t^*]} (u) \le    \int   u(r') P_{r';0}\Big[B_{t^*} \le \eps\Big] \,dr'
 + \int_0^{t^*} jP_{0;s}\Big[ B_{t^*} \le \eps\Big] \,ds
\end{eqnarray*}
Since $B_{t^*}$ has  the law of reflected Brownian motion
 \begin{eqnarray*}
&&
P_{r';0}\Big[B_{t^*} \le \eps\Big]  \le \frac {2\eps}{\sqrt{2\pi t^*}},\;\;\qquad P_{0;s}\Big[ B_{t^*} \le \eps\Big]
\le  \frac{ 2\eps}{\sqrt{2\pi (t^*-s)}}
\end{eqnarray*}
which yields
\begin{eqnarray*}
0&\le & jt^*+F(0;u)- \Delta^{X^{ V}}_{[0,t^*]} (u) 
\le    F(0;u) \cdot \frac {  2\eps}{\sqrt{2\pi t^*}} + \frac{  4j \eps \sqrt{t^*}}{\sqrt {2\pi}}
\end{eqnarray*}
Hence  $\Delta^{X^V}_{[0,t^*]} (u) \to  jt^*+F(0;u)$ as $V\to V^*$.

\medskip We shall next prove that  $\Delta^{X^{ V}}_{[0,t^*]}(u)$   goes to $0$ as $V\to \infty$.
Let $\eps>0$ small,  $V = \eps^{-\frac 34}$ and $r'< X_0 -\eps^{\frac 14}$.  Call $r_k = X_0+
Vt_k$, $t_k= k\eps$, then
 \begin{eqnarray*}
&& P_{r';0}\Big[\tau_0^{X^{ V}} \le  t^*\Big] \le \sum_{k=1}^\infty P_{r';0}\Big[\max_{t\le t_k} B_t\ge r_{k-1}\Big]
\end{eqnarray*}
 Denoting by $P^0_{r';0}$ the law of the Brownian motion on the whole $\mathbb R$ (i.e.\ without reflections at $0$), we have
 \[
 P_{r';0}\Big[\max_{t\le t_k} B_t\ge r_{k-1}\Big] \le 2  P^0_{r';0}\Big[\max_{t\le t_k} B_t\ge r_{k-1}\Big]
 \]
 We then bound
  \begin{equation}
\label{nn5.3.1}
P_{r';0}\Big[\max_{t\le t_k} B_t\ge r_{k-1}\Big] \le 4 P^0_{r';0}\Big[ B_{t_k}\ge r_{k-1}\Big]
 \le {4 e^{-\frac{k}{4 \sqrt \eps}} \int_{k\eps^{\frac 14}}^\infty \frac {e^{-\frac {x^2}{4k \eps}}}{\sqrt{2\pi k \eps}} \le 4\sqrt 2 e^{-\frac{k}{4 \sqrt \eps}}}
%
%
\end{equation}

so that the first term on the right hand side of \eqref{5.3.1} is bounded by
  \begin{eqnarray*}
&&  {\|u\|_{1}}  \eps^{\frac 14} + \, {\|u\|_{\infty}}4\sqrt 2 \, \sum_{k=1}^\infty  e^{- \frac{{k}}{4\sqrt \eps}}
\end{eqnarray*}
which vanishes as $\eps\to 0$. An analogous argument (which is omitted) applies to  the second term on the right hand side of \eqref{5.3.1}.

\vskip.2cm

 \item[ii)]   We  suppose  $V^*<V<V'$ with  $(V'-V)t^*=:\eps$
and $\eps>0$ small enough. To make notation lighter
we shorthand  $X= \{X_t=X_0+Vt\}$ and $X'=\{X'_t=X_0+V't\}$. 
Then by {\eqref{5.3.1},} \eqref{7.18a} and \eqref{7.19a},
   \begin{eqnarray}
    \label{5.100}
&& \Big|\Delta^X_{[0,t^*]}(u) - \Delta^{X'}_{[0,t^*]}(u)\Big| \le  \int_0^{t^*-\eps} g (s) \, P_{X_s;s} \Big[ \tau_{s}^{X'} >t^*\Big]\,ds + R_\eps \\&&
g(s)= \int dr' u(r')\la^{X}_{r',t^*-\eps}(s) +\int_0^s j\kappa^{X}_{s',t^*-\eps}(s)ds'\nn\\&&
R_\eps:= \int   u(r') P_{r';0}\Big[\tau_0^{X}\in [t^*-\eps,t^*]\Big] \,dr'
 + \int_0^{t^*} jP_{0;s}\Big[ \tau^{X}_s\in [t^*-\eps,t^*]\Big] \,ds\nn
	\end{eqnarray}
We are going to prove that there is a function $o(\eps)$
which vanishes as $\eps\to 0$ so that
   \begin{equation}
    \label{5.101}
\sup_{0 \le s \le t^*-\eps} P_{X_s;s} \Big[ \tau_{s}^{X'} >t^*\Big] \le o(\eps)
	\end{equation}
Fix $s \le t^*-\eps$ and
define $\si_s:= \inf \{t \ge s \, : \: B_t \notin (X_s-\eps^{\frac 34}, X_s+\alpha \eps)\}$, with $\alpha > V'+1$,  then
   \begin{eqnarray}
    \label{5.102}
P_{X_s;s} \Big[ \tau_{s}^{X'} >t^*\Big] &\le&  P_{X_s;s} \Big[ \si_s > s +\eps\Big]+
P_{X_s;s} \Big[ B_{\si_s} < X'_{\si_s}; \si_s \le s +\eps\Big] \nn \\
&\le & P_{X_s;s} \Big[ \si_s > s +\eps\Big]+
P_{X_s;s} \Big[ B_{\si_s} = X_s- \eps^{\frac 34}\Big]
	\end{eqnarray}
because, by the choice of $\alpha$, if $B_{\si_s} = X_s +\alpha \eps$ then
$B_{\si_s} > X'_{\si_s}$ as  one can check that
$ X_s +\alpha \eps > X'_{s+\eps}$.

Since $P_{r;s} \Big[ B_{\si_s} = X_s- \eps^{\frac 34}\Big]$ is a linear function of $r$
which has value 1 at $r = X_s- \eps^{\frac 34}$ and is equal
to $0$ at $r=X_s+\alpha\eps$, it then follows that
   \begin{equation}
    \label{5.103}
P_{X_s;s} \Big[ B_{\si_s} = X_s- \eps^{\frac 34}\Big] \le \alpha \, \eps^{\frac 14}
	\end{equation}
On the other hand, since the probability density of $B_{s+\eps}-X_s$ is  $e^{-x^2/(2\eps)}(2\pi \eps)^{-1/2}$
   \begin{equation}
    \label{5.104}
P_{X_s;s} \Big[ \si_s >s + \eps\Big] \le P_{X_s;s} \Big[|B_{s+\eps} - X_s| \le \eps^{\frac 34}\Big]
\le  \frac{2}{\sqrt{2\pi }} \cdot \eps^{\frac 14}
	\end{equation}
so that \eqref{5.101} is proved.  We then have that the first term on the right hand side of \eqref{5.100}
is bounded by:
\[
o(\eps) \int_0^{t^*} g(s)\, ds \le o(\eps) \cdot (F(0;u)+jt^*)
\]

We shall next  bound the probabilities in $R_\eps$.  Call $Y=X _{t^*-\eps}= X_0 + V(t^*-\eps)$, then
   \begin{equation}
    \label{5.105}
P_{r';0}\Big[\tau_0^{X} \in [t^*-\eps,t^*]\Big] \le P_{r';0}\Big[B_{t^*-\eps} \in [Y-\eps^{\frac 14},Y]\Big]
+ \sup_{r'' \le Y-\eps^{\frac 14}}
 P_{r'';t^*-\eps}\Big[ \max _{t\in [t^*-\eps,t^*] } B_t\ge Y\Big]
	\end{equation}
As before we have
   \begin{equation*}
 P_{r';0}\Big[B_{t^*-\eps} \in [Y-\eps^{\frac 14},Y]\Big] \le \frac{\eps^{\frac 14}}{\sqrt{2\pi (t^*-\eps)}}
 	\end{equation*}
Now suppose $r'' \in [0, Y-\eps^{\frac 1 4}]$, then
    \begin{equation*}
  P_{r'';t^*-\eps}\Big[ \max _{t\in [t^*-\eps,t^*] }B_t\ge  Y\Big] \le
  P_{r'';t^*-\eps}\Big[ \max _{t\in [t^*-\eps,t^*] }(B_t - r'')\ge  \eps^{\frac 14} \Big]
 	\end{equation*}
 By the same argument used in \eqref{nn5.3.1}, the latter is bounded by
     \begin{equation*}
  2
  P_{r';t^*-\eps}\Big[ B_{t^*} - r'\ge  \eps^{\frac 14}\Big] \le 2\int_{ \eps^{\frac 14}}^\infty
  \frac{e^{-\frac {x^2}{2\eps}}} {\sqrt{2\pi \eps}} \; dx
  \le 4\sqrt 2e^{-\frac{1}{4 \sqrt \eps}}
 	\end{equation*}

Analogous bounds are proved for $P_{0;s}\Big[ \tau^{X}_s\in [t^*-\eps,t^*]\Big]$, we omit the details.
We have thus proved that also $R_\eps$ is infinitesimal with $\eps$.

\end{itemize}
\vskip.3cm

\noindent
We have thus proved the  identity in \eqref{5.8888}.  The second statement  (i.e. the inequality in  \eqref{5.8888})
follows because $\Delta^{X^V }_{[0,t]}(u)$ is a non-negative,  non decreasing function of $t$ then its maximum in $t \in [0,t^*]$ is  $\Delta^{X^V }_{[0,t^*]}(u)=jt^*$. \qed

\vskip.5cm

\noindent
\col{\bf Proof of  Theorem \ref{exist}, part (a).}
For each positive integer $ n$ consider
the time grid of mesh $\eps_{n}$ and given
$\rho^{(n)}(r,0)$ use the above lemma to construct
$\rho^{(n)}(r,t)$ in the first interval of the grid, $t\le t^*=\eps_{n}$ observing that
at the final time the total mass is exactly equal to the initial one.
We can then use  again the lemma starting from
$\rho^{( n)}(r,{\eps_{n}})$ to construct
$\rho^{(n)}(r,t)$ in the second interval, $t\in [{\eps_{n}},2\cdot { \eps_{n}}]$, at the end of
which we still have conservation of the total mass.
By iteration $\rho^{(n)}(r,t)$ is then defined for all $t\in [0,T]$,
mass is conserved at all discrete times of the grid and
the mass conservation can be violated only in the interior of the
intervals of the time grid. By Lemma \ref{ciccio}, at any time $t\in [0,T]$
we have conservation of mass modulo at most  $j {\eps_{n}}$.
\qed

\section{Characterization and uniqueness}
\label{sec:5}

{
In this section we characterize the function $\bar \rho$ of Theorem \ref{exist} as the unique separating element between barriers thus proving the existence and uniqueness of the generalized solution.}
%
%
%
%
%
%
%
{
\subsection{The key inequality}
\label{sec:5.1}

We fix a function $\rho_0$, $T>0$ and a quasi-solution $(X_t,\rho(\cdot,t), {\eps})$, $t\in [0,T]$ with accuracy parameter
$\eps$ as in  Definition \ref{defin1.2}. Recalling Definition \ref{def:2.1} we prove the following key inequality.}

%

\vskip.5cm

{
\begin{proposition}
\label{keyes}
For any $\delta>0$, there is a constant $c$ so that  for all $k\in \N$ such that $k \delta \le T$
 \begin{equation}
\label{5.8}
S^{(\delta,-)}_{k\delta} ( \rho(\cdot,0))  \le   \rho(\cdot, k\delta) \le
S^{(\delta,+)}_{k\delta} (  \rho(\cdot,0))  \;\;\;\qquad \text{modulo}\quad ck\eps
  \end{equation}
\end{proposition}

 \noindent {\bf Proof.}
We prove }\eqref{5.8} by induction on $k$.  We thus fix $k$,
shorthand
 \begin{equation}
\label{5.9}
v^{(\pm)} := S^{(\delta,\pm)}_{k\delta} (  \rho(\cdot,0)),\;\qquad\; w=  \rho(\cdot, k\delta),\;\qquad
{Z_t:=}X_{k\delta+t},\qquad
\;m=ck\eps
\end{equation}
and
suppose (by induction) that
 \begin{equation}
\label{5.10}
v^{(-)} \le w \;\;\text{modulo}\; m\qquad \qquad\text{and}  \qquad \qquad \: w \le v^{(+)} \;\;\text{modulo}\; m
\end{equation}
We call  $w (r,t)=\rho(r,k\delta+t)$, $t\in [0,\delta]$.
We shall prove in  Subsection \ref{sec:6a} that
\begin{equation}
\label{5.11}
S^{(\delta,-)}_{\delta}(v^{(-)})\le w(\cdot, \delta)\;\;\;\qquad \text{modulo}\;\; m + c\eps
\end{equation}
and  in Subsection \ref{sec:7a} that
\begin{equation}
\label{5.12}
 w(\cdot, \delta)\le S^{(\delta,+)}_{\delta}(v^{(+)})\;\;\; \qquad \text{modulo}\;\; m + c\eps
\end{equation}
which will thus prove the induction {and concludes the proof of the Proposition. \qed}

\vskip.3cm
\noindent
{We conclude this subsection by observing two properties of $Z_t$ and $w(\cdot,t)$ (see \eqref{5.9})  that, together with  \eqref{1.9}, will be extensively used in next two subsections for the proofs of the upper and lower bound.}

\vskip.3cm
\noindent
\colo{{$(Z_t,w(\cdot,t),2\eps)$}, $t\in [0,\delta]$,  is a
quasi-solution of the FBP \ref{defin1.1} with initial datum $w$
and accuracy parameter $2\eps$ (see \eqref{1.9}), then  from \eqref{3.144444} and \eqref{3.1444441111111} it follows that for any $I \subseteq \mathbb R^+$,
 \begin{eqnarray}
\label{5.6}
&& \hskip-1cm \int_I w(r,t)\,dr =   \int w(r') P_{r';0}\Big[\tau_0^Z >t; B_t\in I\Big] \,dr'
+ \int_0^t jP_{0;s}\Big[ \tau^Z_s >t; B_t\in I\Big] \,ds
\end{eqnarray}
Moreover if $\Delta^Z_{[0,t]}(w)$   is the  mass lost by $w$  in the time interval $[0,t]$, $0\le t \le \delta$ (see  \eqref{5.3.1}), then from \eqref{5.3.1.1} and \eqref{1.9} we have
\begin{equation}
\label{5.7}
 {\sup_{t\in [0,\delta]} \Big| \Delta^Z_{[0,t]}(w)- jt\Big| \le 2\eps}\;.
\end{equation}}

\subsection{Lower bound (proof of  \eqref{5.11}).}
\label{sec:6a}

We write here $v$ for $v^{(-)}$ and by the induction hypothesis we have that $v\le w$ modulo $m$.
We need to prove that for all $r\ge 0$
	\begin{equation}
	\label{n5.10a}
 F(r;K^{(\delta)}G_\delta^{\rm neum}*v)\le F(r;w(\cdot,\delta))+m+c\eps
\end{equation}
We decompose $v=v^*+v_1$ with
\[
v^*= v 	\,{\mathbf 1_{[0, R_m]}},\;\qquad \quad  R_m: \int_{R_m}^\infty v =m,\qquad \quad v_1 = v \, {\mathbf 1_{(R_m,+\infty]}}
\]
so that, by
{Lemma \ref{prop4.2}}, $v^*\le w$.  Let
\[
w=\tilde w+w_1,\;\;\quad w_1= w {\mathbf 1_{[0,\tilde R]}},\quad\qquad \tilde R: \int w_1 = F(0;w)-F(0;v^*)\ge 0
\]
(because $v^*\le w$),
so that by {Lemma \ref{lemma4.2},}
\[
v^*\le \tilde w,\qquad F(0;v^*)=F(0;\tilde w)
\]
{From \colo{the right identity in }\eqref{3.1} we have }
 \begin{equation}\label{Fv}
 F(r;G_\delta^{\rm neum}*v^*)=\int   dr' v^*(r')
  \int_r^\infty
  G_\delta^{\rm neum}(r', z) \, dz=
  \int   \;   v^*(r') \; P_{r'{;0}}\Big[ B_\delta \ge r\Big]dr'
 \end{equation}
\colo{We define   $\tilde w (r,t)$, $t\in [0,\delta]$, the evolution at time $t$ of the profile $\tilde w(r)$ with  the dynamics defined by \eqref{5.6},}
 then
%
 \begin{eqnarray}
\label{6.2}
F(r;\tilde w(\cdot,\delta)) =
\int \; \tilde w(r') \, P_{r'{;0}}\Big[{\tau_0^Z>\delta ;}\; B_\delta \ge r\Big] dr' + \int_0^\delta j  \, P_{0;s}\Big[{\tau_s^Z>\delta ;}  \;B_\delta \ge r\Big]ds
\end{eqnarray}
Let $\Delta^Z_{[0,\delta]}(w)$ and $\Delta^Z_{[0,\delta]}(\tilde w)$  be the  mass lost by $w$ and $\tilde w$ in the time interval $[0,\delta]$. Since {$\tilde w(r) \le w(r)$ for all $r$,
then, from \eqref{5.3.1} it follows that $\Delta^Z_{[0,\delta]}(\tilde w) \le \Delta^Z_{[0,\delta]}(w)$, then, from \eqref{5.7} we have
\begin{equation}
\label{6.3}
\colo{{\Delta^Z_{[0,\delta]}(\tilde w)}  \le {\Delta^Z_{[0,\delta]}(w)} \le  j\delta +{2\eps}}
\end{equation}
{From \eqref{5.3.1} and  \eqref{6.2} we get}
 \begin{eqnarray}
\label{6.4}
F(r;\tilde w(\cdot,\delta)) &\ge &
\int \; \tilde w(r') \, P_{r'{;0}}\Big[ B_\delta \ge r\Big]dr'
+\int_0^\delta j  P_{0;s}\Big[ B_\delta \ge r\Big] ds- {\Delta^Z_{[0,\delta]}(\tilde w)}
\end{eqnarray}
{The function  $ \phi(r'):=P_{r'{;0}}\Big[ B_\delta \ge r\Big]$ is bounded and non decreasing, then from \eqref{4.7aaa}
\begin{equation}
\int \; v^*(r') \, \phi(r')dr' \le \int  \; \tilde w(r') \, \phi(r') dr'
\end{equation}}
hence
 \begin{eqnarray*}
F(r;\tilde w(\cdot,\delta)) \ge
\int  v^*(r') P_{r'{;0}}\Big[ B_\delta \ge r\Big]dr'
+\int_0^\delta j \,  P_{0;s}\Big[ B_\delta \ge r\Big] ds- {\Delta^Z_{[0,\delta]}(\tilde w)} &&\nn \\
{\ge
F(r;G_\delta^{\rm neum}*v^*)
+\int_0^\delta j  \, P_{0;s}\Big[ B_\delta \ge r\Big] ds- (j\delta +  \col{2}\eps)} &&
\end{eqnarray*}
where the last inequality follows from \eqref{Fv} and \eqref{6.3}. Then
 \begin{eqnarray*}
&& \tilde w(\cdot,\delta)\ge G_\delta^{\rm neum}*v^*\;\;\; \qquad\text{modulo}\;\; j\delta +   \col{2}\eps   \\ \text{and} \qquad &&
F(0;\tilde w(\cdot,\delta)) \ge F(0;G_\delta^{\rm neum}*v^*) - \col{2}\eps \nn
\end{eqnarray*}
{Using \eqref{4.15a} with $m=j\delta+\col{2}\eps$ and $\alpha =\col{2}\eps$ we get}
 \begin{equation*}
\tilde w(\cdot,\delta)\ge K^{(\delta)}G_\delta^{\rm neum}*v^*\;\;\; \qquad \text{modulo}\;\; \col{2}\eps
\end{equation*}
On the other hand
 \begin{equation*}
 K^{(\delta)}G_\delta^{\rm neum}*v^*\ge K^{(\delta)}G_\delta^{\rm neum}*v\;\;\; \qquad \text{modulo}\;\;
 {F(0;K^{(\delta)}G_\delta^{\rm neum}*v_1)}
\end{equation*}
where {$F(0;K^{(\delta)}G_\delta^{\rm neum}*v_1) =F(0; v_1)=m$} so that
  \begin{equation*}
\tilde w(\cdot,\delta)\ge K^{(\delta)}G_\delta^{\rm neum}*v\;\;\;\qquad \text{modulo}\;\; m+  \col{2}\eps
\end{equation*}
which proves \eqref{n5.10a}
 (with $c\ge 1$)
because
$w(\cdot,\delta) \ge \tilde w(\cdot,\delta)$.
\qed

\vskip1cm

\subsection{Upper bound (proof of  \eqref{5.12}).}
\label{sec:7a}

We shorthand here $v$ for $v^{(+)}$, \colo{then, by the induction assumption, $w \le v$ modulo $m$ where $w(r,t)$, $t\in [0,\delta]$, is  given by \eqref{5.6}. We have to prove that
\begin{equation}
\label{7.1a}
I(r):= F(r; w(\cdot, \delta)) - F(r; G_\delta^{\rm neum}*K^{(\delta)}v)- m  \le c\eps,\quad \text{for all $r\ge 0$}
\end{equation}}
\noindent
We write $w=w_0+w_1+w_2$ with
	\begin{equation}
	\label{7.2a}
w_2= w \, {\mathbf 1_{(R_m(w),+\infty)}}:\;\;\;\int w_2=m;
\qquad \qquad w_1=w\,\mathbf 1_{[R_1,R_m(w)]}:\;\; \int w_1 = j\delta
		\end{equation}
Next lemma shows that after some simple manipulations the proof of the upper bound
is reduced to some inequalities among the $w_i$.

\medskip

\begin{lemma}
 Let $I(r)$ be the function defined in \eqref{7.1a}, then 
\begin{eqnarray}
&&  I(r)\le \int w_1(r')P_{r';0}\Big[B_\delta \ge r; {\tau_0^Z} >\delta\Big]dr' -
\sum_{i\in\{0,2\}} \int w_i(r') P_{r';0}\Big[B_\delta \ge r; {\tau_0^Z}  \le \delta\Big]dr' \nn
\\&& \hskip1cm - j\int_0^\delta  P_{0;s} \Big[B_\delta \ge r; {\tau_s^Z}  \le \delta\Big] ds
\label{7.12a}
	\end{eqnarray}

\end{lemma}

\medskip

\noindent
{\bf Proof.}
Since $w\le v$ modulo $m$, by \eqref{4.10a} (and the definition of $w_2$)
$w_0+w_1 \le v$.
\colo{We  write  $v=v_0+v_1+ v_2$ where
	\begin{equation}
	\label{7.3abbb}
v_1=v\, {\mathbf 1_{[R_\delta(v), +\infty)}}:\;\; \int v_1 = j\delta
		\end{equation}
\[
v_2= v \, {\mathbf 1_{[0,\tilde R]}}:\quad\;\;\int v_2=F(0;v)-F(0;w_0+w_1)
\]
so that by \eqref{4.8a}
\[
w_0+w_1 \le v_0+v_1,\quad F(0; w_0+w_1)=F(0; v_0+v_1)
\]
then $$F(0;v_i) = F(0;w_i),\qquad i=0,1$$
 hence, by the last inequality in \eqref{4.12a},
	\begin{equation}
	\label{7.3a}
w_0 \le v_0,\qquad F(0;v_0) = F(0;w_0)
		\end{equation}
Observe moreover that
	\begin{equation}
	\label{7.3abc}
F(r;S^{(\delta,+)}_{\delta}(v_0+v_1)) \le
F(r;S^{(\delta,+)}_{\delta}(v)), \qquad   {S^{(\delta,+)}_{\delta}(v_0+v_1)= G^{\rm neum}_\delta *[v_0
+ j \delta D_0]}
		\end{equation}}
We call
\begin{eqnarray}
\label{7.6a}
&& f_i (r)= \int \; w_i(r')\, P_{r';0}\Big[B_\delta \ge r; \; {\tau_0^Z} >\delta\Big]dr' ,\;\qquad i=0,1,2\\&&
f_3(r)=	j\int_0^\delta  \; P_{0;s} \Big[B_\delta \ge r;\; {\tau_s^Z}  >\delta\Big]	ds
\nn
	\end{eqnarray}
so that, \colo{from \eqref{5.6} we have}
	\begin{equation}
	\label{n4.21}
F(r; w(\cdot, \delta)) = \sum_{i=0}^3{f_i(r)}
	\end{equation}
For $i\in\{0,1,2\}$ we write
\begin{eqnarray}
&& f_i (r)= \int w_i(r')P_{r';0}\Big[B_\delta \ge r\Big] dr' -\int  w_i(r') P_{r';0}\Big[B_\delta \ge r; {\tau_0^Z}  <\delta\Big]dr'
\label{7.6aaa}
	\end{eqnarray}
and  for $i=3$:
\begin{eqnarray}
\label{7.7a}
&& 
{f_3(r)=	j}\int_0^\delta  P_{0;s} \Big[B_\delta \ge r\Big] ds-j
\int_0^\delta P_{0;s} \Big[B_\delta \ge r; \; {\tau_s^Z}  <\delta\Big]ds
\label{7.8a}
	\end{eqnarray}
Using \eqref{n4.21} and the left inequality in \eqref{7.3abc} we then get
	\begin{equation}\label{S1}
I(r)\le \sum_{i=0}^3f_i(r)-F(r; S^{(\delta,+)}_\delta(v_0+v_1))-m
	\end{equation}	
with $f_i (r)$ as in \eqref{7.6aaa}--\eqref{7.8a} when $i\ne 1$ and by $f_1(r)$
	as in \eqref{7.6a}.
From the  identity in \eqref{7.3a} and \eqref{3.1} we get
\begin{equation}
\label{7.9a}
F(r; S^{(\delta,+)}_\delta(v_0+v_1))=F\Big(r;G^{\rm neum}_{\delta} *[v_0(\cdot,\delta)+j\delta D_0]\Big) =
\int v_0(r')P_{r';0}\Big[B_\delta \ge r\Big] dr' +j\delta P_{0;0}\Big[B_\delta \ge r\Big]
\end{equation}
\colo{From \eqref{4.11a} and \eqref{3.1} we have}
\begin{eqnarray}
\label{7.10a}
&&  \int w_0(r')P_{r';0}\Big[B_\delta \ge r\Big]dr'  \le  \int  v_0(r')P_{r';0}\Big[B_\delta \ge r\Big]dr'
\end{eqnarray}
moreover
\begin{eqnarray}\label{S2}
&&	j\int_0^\delta  P_{0;s} \Big[B_\delta \ge r\Big] ds \le j\delta P_{0;0}\Big[B_\delta \ge r\Big] \nn \\ \text{and} \quad &&  \colo{\int w_2(r')P_{r';0}\Big[B_\delta \ge r\Big]dr'  \le  m}
	\end{eqnarray}
 \colo{then \eqref{7.12a} follows from \eqref{S1}, \eqref{7.9a}, \eqref{7.10a} and \eqref{S2}.}   \qed

\vskip.5cm

%
%
%
%
%
%
%
%
%
%
%
%
\noindent
From the above lemma we are left with the proof that
\begin{eqnarray}
&&  \int w_1(r')P_{r';0}\Big[B_\delta \ge r; \; {\tau_0^Z}  >\delta\Big]dr'  \le
\sum_{i\in\{0,2\}} \int  w_i(r') P_{r';0}\Big[B_\delta \ge r; \;{\tau_0^Z}  \le \delta\Big]dr' \nn
\\&& \hskip1cm + j\int_0^\delta \; P_{0;s} \Big[B_\delta \ge r; {\tau_s^Z}  \le \delta\Big] ds + c\eps
\label{7.12abb}
	\end{eqnarray}
Call
\begin{equation}
\label{7.13a}
\alpha(r):= P_{r;0}\Big[{\tau_0^Z}  \le \delta\Big],\qquad \beta(s):= P_{0;s}\Big[{\tau_s^Z}  \le \delta\Big]
\end{equation}
Then the inequality in \eqref{7.12abb} becomes
\begin{eqnarray}
&& \hskip-1cm \int w_1(r')[1-\alpha(r')]P_{r';0}\Big[B_\delta \ge r\,\big|\,
{\tau_0^Z} >\delta\Big]dr' \le
\sum_{i=0,2} \int  w_i(r') \alpha(r')P_{r';0}\Big[B_\delta \ge r\,\big|\, {\tau_0^Z} \le \delta\Big]dr' \nn
\\&& \hskip1cm + j\int_0^\delta  \beta(s) P_{0;s}\Big[B_\delta \ge r\,\big|\, {\tau_s^Z}  \le\delta\Big] ds+ c\eps
\label{7.14a}
	\end{eqnarray}
	By using \eqref{5.3.1.1} and \eqref{5.7}, we get
		\begin{equation}
		\label{n5.28a}
\colo{\big|F(0, w(\cdot,\delta))-F(0,w)\big|=\big|\Delta^Z_{[0,\delta]} (w)-j\delta\big|
\le 2\eps}
		\end{equation}
\colo{From \eqref{5.6} we have}
	\begin{equation*}
	\colo{F(0,w(\cdot,\delta))-F(0,w)}=j\delta-
	\sum_{i=0}^2\int \alpha(r)w_i(r) dr\col{-} j\int_0^\delta \beta(s) ds
	\end{equation*}
Since  $\int w_1 = j\delta$ we have
\begin{equation}
\label{7.15a}
\Big| \Big\{\sum_{i=0}^2\int \alpha(r)w_i(r) dr+ j\int_0^\delta \beta(s) ds\Big\} - \int w_1(r) dr\Big| \le 2\eps
\end{equation}
\colo{Since \eqref{7.14a} is trivially satisfied for any $c\ge 2$ when} $\int w_1(r)[1-\alpha(r)] \, dr  \le  2\eps$, we can  suppose that $\int w_1(r)[1-\alpha(r)] \, dr >  2\eps$. It follows that  there is an interval $I\in \mathbb R_+$ so that
\[
\int_I w_1(r)[1-\alpha(r)] \, dr  = 2\eps
\]
By \eqref{7.15a} there exists $q \le 1$ so that}
\begin{equation}
\label{7.15aaa}
M:=\int_{I^c} w_1(r)[1-\alpha(r)] =q \Big\{\sum_{i=\col{0,2}}\int \alpha(r)w_i(r) + j\int_0^\delta \beta(s) \Big\}
\end{equation}
We are going to prove that there exists a constant $c'$ so that
			\begin{eqnarray}
			\nn
&&\hskip-1cm\int_{I^c}  w_1(r')[1-\alpha(r')]P_{r';0}\Big[B_\delta \ge r\,\big|\, {\tau_0^Z}  >\delta\Big]dr'\le
q\Big(\sum_{i=0,2} \int  w_i(r') \alpha(r')P_{r';0}\Big[B_\delta \ge r\,\big|\, {\tau_0^Z}  \le \delta\Big] dr'\\&&\hskip6cm
  + j\int_0^\delta  \beta(s) P_{0;s}\Big[B_\delta \ge r\,\big|\, {\tau_s^Z}  \le\delta\Big] \Big)ds+ c'\eps\nn\\
\label{7.14abcd}
	\end{eqnarray}
which yields \eqref{7.14a} with $c=c'+2$.

\noindent
Recalling \eqref{7.18a}--\eqref{7.19a} there exists a non negative measure $g(dt)$ on $[0,\delta]$,
so that $\dis{\int_0^\delta g(dt) =M}$ and
\begin{eqnarray*}
&& \hskip-1cm q\Big(\sum_{i=0,2} \int dr' w_i(r') \alpha(r')P_{r';0}\Big[B_\delta \ge r\,\big|\, {\tau_0^Z}  \le \delta\Big] \\&&
\hskip1cm + j\int_0^\delta ds \beta(s) P_{0;s}\Big[B_\delta \ge r\,\big|\, {\tau_s^Z}  \le\delta\Big] \Big)=
\int_0^\delta  \, g(dt) P_{ Z_t;t}\Big[B_\delta \ge r\Big] 
	\end{eqnarray*}

\noindent
Since $w_1(r)[1-\alpha(r)]\mathbf 1_{I^c}(r)dr=:\mu(dr)$ and $g(dt)$ have same mass $M$,
and $g(dt)$  does not have atomic components, then,
by the isomorphism of Lebesgue measures, \cite{roklin},  there is a map $\Ga:\mathbb R_+ \to [0,\delta]$ so that
\begin{eqnarray}
&&
\int_0^\delta g(dt) P_{ Z_t;t}\Big[B_\delta \ge r\Big]  = \int \mu(dr')
P_{ Z_{\Ga(r')};\Ga(r')}\Big[B_\delta \ge r\Big]
\label{7.24a}
	\end{eqnarray}
\eqref{7.14abcd} with $c'=0$ will then follow from the inequality
\begin{eqnarray}
&&
 P_{r';0}\Big[B_\delta \ge r\,\big|\, {\tau_0^Z} >\delta\Big] \le
P_{Z_{t};t}\Big[B_\delta \ge r\Big],\quad r' \in [0, Z_0),\; t\in [0,\delta)
\label{7.25a}
	\end{eqnarray}
(used with  $t= \Ga(r')$).  A proof of
\eqref{7.25a} via a coupling  between conditioned and unconditioned Brownian
motions which preserves
order is given in \cite{CDP}.  Here we present a more elementary proof where we use
classical couplings between Brownian motions (reflected at the origin).

\begin{proposition}
\label{prop5.1}
Let $P_{r,r',s} = P_{r,s}\times P_{r',s}$, $r,r'\in \mathbb R_+$,  be the product of two independent
standard Brownians, $B_t^{(1)}$
and $B_t^{(2)}$, $t\ge s$, on $\mathbb R_+$ with reflections at $0$.  Call $\tau:= \inf\{s: B_s^{(1)}=
B_s^{(2)}\}$ (and equal to $+\infty$ if
the set is empty).  Then
\begin{equation}
\label{5.4}
b_t := \begin{cases} B_t^{(2)} & \text{if}\; t\le \tau \\
 B_t^{(1)} & \text{if}\; t\ge \tau \end{cases}
\end{equation}
has the law of a Brownian motion starting from $r'$ at time $s$. 

\end{proposition}

\medskip

\noindent
{\bf Proof of \eqref{7.25a}.} We use the method of duplicating variables.
Let $\ga^{-1}$ be a positive integer (eventually $\ga\to 0$), $B^{(1)}_i$, $i=1,..,\ga^{-1}$
independent Brownian motions which start moving at time $t$  from $Z_t$ (and are frozen
before: $B^{(1)}_i(s)= Z_t$, $s\le t$).  Then
denoting below by $E^{(1)}$ the expectation with respect to such independent Brownians,
\begin{eqnarray}
&& P_{ Z_{t};t}\Big[B_\delta \ge r\Big]=
E^{(1)}\Big[ \ga \sum_i {\mathbf 1_{[r, +\infty)}(B^{(1)}_i(\delta))}\Big]
\label{7.26a}
	\end{eqnarray}
We can proceed in analogous way with $ P_{r';0}\Big[B_\delta \ge r\,\big|\,{\tau_0^Z}>\delta\Big]$
which is now conveniently rewritten as
\begin{eqnarray}
&& P_{r';0}\Big[B_\delta \ge r\,\big|\, {\tau_0^Z} >\delta\Big]=  P_{r';0}\Big[B_\delta \ge r\,;\,{\tau_0^Z}  >\delta\Big]
\Big( 1 + \Big\{ \frac{1}{1-\alpha(r')} - 1\Big\} \Big)
\label{7.27a}
	\end{eqnarray}
Calling $N_\ga := $ the integer part of $\ga^{-1}\big\{ \frac{1}{1-\alpha(r')} - 1\big\}$,
we then consider  $B^{(2)}_i$, $i=1,..,\ga^{-1}+N_\ga$ 
independent Brownian motions which start  at time $0$  from $r'$ and are removed once they reach
the edge $Z_t$.  Then, denoting below by $E^{(2)}$ expectation with respect to the law of the independent Brownians
with death at the edge,
\begin{eqnarray}
&& P_{r';0}\Big[B_\delta \ge r\,\big|\, {\tau_0^Z} >\delta\Big]= \lim_{\ga\to 0}E^{(2)}\Big[ \ga \sum_{i=1}^{\ga^{-1}+N_\ga} {\mathbf 1_{[r,+\infty)}(B^{(2)}_i(\delta))}\Big] 
\label{7.28a}
	\end{eqnarray}
The equality holds only in the limit because of the integer part in the definition of $N_\ga$.

Call $i_1,..,i_{M_\ga}$ the labels of the $B^{(2)}$-particles still existing at time $t$. To simplify notation
we relabel them as $1,..,M_\ga$ and denote by $r_1,..,r_{M_\ga}$ the  positions of the corresponding
$B^{(2)}$-particles at time $t$.  Thus at time $t$ we have $\ga^{-1}$  $B^{(1)}$-particles all at
$Z_t$ and a number $M_\ga$ of $B^{(2)}$-particles at $r_1$,..,$r_{M_\ga}$ (all $<  Z_t$,
otherwise they would be dead).

We call married couples at time $t$ the pairs  $(B^{(1)}_i(t),B^{(2)}_i(t))$ with
$i\le \min\{\ga^{-1},M_\ga\}$, single the $B^{(2)}_i(t)$ particles with $i>\ga^{-1}$, if  any.
We then let evolve all the particles present at time $t$ in the following way: each married
couple moves as in Proposition \ref{prop5.1} independently of the other couples and
of all the other particles which move as independent Brownians.  At the first time $s>t$ when a $B^{(2)}_i$
particle reaches the edge $Z_s$, it disappears.  If it was single then the evolution after $s$ proceeds
as before without the dead particle.  If instead the particle which reaches the edge at time $s>t$ is  one of the ``married particles'' the $B^{(1)}$-particle in the pair is coupled with
a single $B^{(2)}$-particle, if such a particle exists at time $s$, otherwise
it is an unmarried  $B^{(1)}$-particle.  After that the process continues with same rules
and it is thus  defined completely by iteration.
We denote by $P$ and $E$ law and expectation of the coupled process.

By its definition   it  follows that the number $K_\ga$ of single {$B^{(2)}$-particles}
at the final time $\delta$ is
\[
K_\ga = \max\Big\{0;N_\ga - \sum_{i=1}^{\ga^{-1}+N_\ga} \mathbf 1_{{B^{(2)}_i (s)=Z_s, \;\text{for some $s\in [0,\delta]$}}}\Big\}
\]
while for all pairs $(B_i^{(1)}(\delta),B_i^{(2)}(\delta))$ present at time $\delta$ the inequalities
$B_i^{(1)}(\delta)\ge B_i^{(2)}(\delta)$ hold.  Thus by \eqref{7.26a}--\eqref{7.28a}
\begin{eqnarray}
&& P_{ Z_{t};t}\Big[B_\delta \ge r\Big]- P_{r';0}\Big[B_\delta \ge r\,\big|\,{\tau_0^Z}>\delta\Big] \ge
{-}\lim_{\ga\to 0} E\Big[\ga K_\ga \Big]
\label{7.29a}
	\end{eqnarray}
By the law of large numbers for independent variables,
for any $\zeta>0$
\begin{equation}
\lim_{\ga \to 0} P_{r',0}\Big[\; \Big| \ga\sum_{i=1}^{\ga^{-1}+N_\ga} \mathbf 1_{{B^{(2)}_i (s)=Z_s, \;\text{for some $s\in [0,\delta]$}}}
- (1+\ga N_\ga) P_{r',0}[ {\tau_0^Z} \le \delta] \Big| \le \zeta\Big] = 1
\label{7.30a}
	\end{equation}
Since
\[
\lim_{\ga \to 0} (1+\ga N_\ga) P_{r',0}[ {\tau_0^Z} \le \delta] = \frac{\alpha(r')}{1
-\alpha(r')} = \lim_{\ga \to 0} \ga N_\ga
\]
it follows from \eqref{7.30a} that
\begin{equation*}
\lim_{\ga\to 0} P_{r',0}\Big[ \;\Big| \ga\sum_{i=1}^{\ga^{-1}+N_\ga} \mathbf 1_{{B^{(2)}_i (s)=Z_s, \;\text{for some $s\in [0,\delta]$}}} -\ga N_\ga\Big|\le \zeta\Big] =1
\label{7.31a}
	\end{equation*}
hence $\dis{\lim_{\ga\to 0} P\big[\ga K_\ga \le \zeta \big]=1}$ which yields
$\dis{\lim_{\ga\to 0} E\big[\ga K_\ga  \big]=0}$, thus  the right hand side of
\eqref{7.29a} is equal to 0.

\qed

\vskip1cm

\subsection{Conclusion of the proof of Theorem \ref{exist} and  Theorem \ref{ineq}}
\label{sec:8a}

{
In all this subsection  $(X^{(n)}_t,\rho^{(n)},\eps_n)$, $t\in (0,T)$  is an optimal sequence of quasi-solutions
as in part (b) of  Theorem \ref{exist}. The following result is an easy consequence of Proposition \ref{keyes}.
	\begin{corollary}
	\label{coro:5.4}
For any fixed $t\in (0,T)$ the sequence of non negative measures $\{\rho^{(n)}(r,t) dr\}$ on $\mathbb R_+$
with its Borel $\si$-algebra
is
tight
(in the weak topology of measures).
\end{corollary}

\noindent {\bf Proof.}
Fix arbitrarily
$t\in (0,T)$ and $\ell\in \mathbb N$, call  $\delta = 2^{-\ell} t$ and $k^*= 2^\ell$ so that $k^*\delta =t$. Then from \eqref{5.8} with such value of $\delta$ and with $k=k^*$, we get that
 there exists a constant $c$ so that for any $n$
 \begin{equation}
\label{n4.1}
S^{(\delta,-)}_{t} (  \rho^{(n)}(\cdot,0))  \le   \rho^{(n)}(\cdot, t) \le
S^{(\delta,+)}_{t} (  \rho^{(n)}(\cdot,0))  \;\;\; \qquad \text{modulo}\quad {c  k^*\eps_n}
\end{equation}
%
To prove tightness we use the
Prokhorov theorem: it is then sufficient to show that
for any $\zeta>0$ there is $r_\zeta$ so that
  \begin{equation}
F(r_\zeta;\rho^{(n)}(\cdot,t)) \le \zeta\quad \text{for all $n$}
\label{8.1a}
	\end{equation}
	From \eqref{n4.1} it follows that}
  \begin{equation}
F(r;\rho^{(n)}(\cdot,t)) \le  F(r;S^{(\delta,+)}_t(\rho^{(n)}(\cdot,0))) + c{\color{brown}k^*} \eps_{n}
\label{8.2a}
	\end{equation}
By \eqref{3.10} {and the definition of quasi-solution
  \begin{eqnarray}
&&\hskip-1cm \big| F(r;S^{(\delta,+)}_t(\rho^{(n)}(\cdot,0))) - F(r;S^{(\delta,+)}_t(\rho_0))\big|  \le
\int \big|S^{(\delta,+)}_t(\rho^{(n)}(\cdot,0))-S^{(\delta,\pm)}_t(\rho_0)\big|(r) \, dr \nn\\&&
\hskip2cm\le
\int
 \big|\rho^{(n)}(r,0)-\rho_0(r)\big|\, dr\le \eps_n
\label{8.3a}
	\end{eqnarray}}
%
%
Moreover {for any $\zeta$} there is $R$ so large that for all $r\ge R$
  \begin{equation}
F(r;S^{(\delta,+)}_t(\rho_0))) \le 2jt\int_r^\infty \frac{e^{-x^2/2t}}{\sqrt{2\pi t}}dx
 +\int_r^\infty dx \int_{\mathbb R_+} dr' \rho_0(r') G^{\rm neum}(r',x)
 \le \frac \zeta 2
\label{8.6a}
	\end{equation}
{
Given  $\zeta$ let ${n_0}$ be  so  that for all $n\ge  {n_0}$, $\eps_{n}(1+c\delta^{-1}T)\le \zeta/2$, then \eqref{8.1a} follows from \eqref{8.2a}, \eqref{8.3a} and \eqref{8.6a} and thus the corollary is proved.\qed
\vskip.4cm

\noindent
{\bf Proof of (b) and (c) of Theorem \ref{exist}}. From Corollary \ref{coro:5.4}
the sequence  $\{\rho^{(n)}(r,t) dr\}_{n\ge 0}$} converges weakly by subsequences to a limit
measure $\mu(dr)$ and we next prove that $\mu(dr)= S_t(\rho_0) dr$. {This  implies that
the sequence itself converges and it identifies the function $\bar \rho$ in (c), thus concluding both the proof of Theorem \ref{exist} and the first statement of Theorem \ref{ineq}. }

Call
  \begin{equation}
f'_{r,t} : = \liminf_{n\to \infty} F(r;\rho^{(n)}(\cdot,t)),\quad
 f''_{r,t} : = \limsup_{n\to \infty} F(r;\rho^{(n)}(\cdot,t))
 \label{8.7a}
	\end{equation}
Let again $\delta= 2^{-\ell}t$, $\ell\in \mathbb N$,
then {by \eqref{n4.1}}
  \begin{equation*}
 F(r;S^{(\delta,-)}_t(\rho^{(n)}(\cdot,0))) -c\eps_{n} \delta^{-1}T \le F(r;\rho^{(n)}(\cdot,t)) \le  F(r;S^{(\delta,+)}_t(\rho^{(n)}(\cdot,0))) + c \eps_{n} \delta^{-1}T
	\end{equation*}
By \eqref{8.3a} 
  \begin{equation}
 F(r;S^{(\delta,-)}_t(\rho_0)) -\zeta_{n}  \le F(r;\rho^{(n)}(\cdot,t)) \le  F(r;S^{(\delta,+)}_t(\rho_0)) +  \zeta_{n} \label{8.8a}
	\end{equation}
where {$\zeta_{n}:= \eps_{n}(c \delta^{-1}T + 1)$}.
Since $\zeta_{n}\to 0$ as $n\to \infty$
and recalling that $\delta=2^{-\ell}t$:
  \begin{equation}
 F(r;S^{(2^{-\ell}t,-)}_t(\rho_0))  \le f'_{r,t}  \le f''_{r,t} \le  F(r;S^{(2^{-\ell}t,+)}_t(\rho_0))
 \label{8.9a}
	\end{equation}
By \eqref{4.4a}  letting $\ell\to \infty$
  \begin{equation}
  f'_{r,t}  = f''_{r,t} =  F(r;S_t(\rho_0))
 \label{8.10a}
	\end{equation}
which, by the arbitrariness of $r$,  identifies the limit measure $\mu(dr) = S_t(\rho_0)(r)dr$. \qed

\vskip.4cm

{
\noindent {\bf Proof of Theorem \ref{ineq}}. The first statement has been already proved above. }
The function  $\rho(r,t) :=\rho_a(r)$, see \eqref{1.6}, is a classical (and stationary) solution
in the sense of Definition \ref{defin1.1}, hence it is a fortiori a quasi-solution (with $X_t \equiv \frac{a}{2j}$ and accuracy parameter $\eps=0$). Then for what proved so far, $S_t(\rho_a)=\rho_a$.
Let $\rho_0$ be
as in Theorem \ref{ineq}, then, since $R(\rho_0)<\infty$ for $a$ large enough $\rho_0(r)\le \rho_a(r)$ for all $r$
hence $\rho_0\le \rho_a$ also in the sense of mass transport.  Since $S_t$ preserves order (see the last statement in Theorem \ref{thmineq}), we have
$S_t(\rho_0) \le S_t(\rho_a)=\rho_a$ which proves that $R(\rho(\cdot,t) )\le R(\rho_a)=\frac{a}{2j}{:=a_2}$ for all $t\ge 0$.
{Analogously, since $R(\rho_0)>0$ for $\bar a$ such that $\frac{\bar a}{2j}<R(\rho_0)$ we get that $\rho_{\bar a}=S_t(\rho_{\bar a})\le S_t(\rho_0)$
which proves that $R(\rho(\cdot,t) )\ge R(\rho_{\bar a})=\frac{\bar a}{2j} :=a_1$ for all $t\ge 0$. }
\qed

\vspace{0.5cm}

{\bf Acknowledgments.} The authors thank F. Comets, S. Luckhaus, J. R. Ockendon, D. Ioffe, S. Olla, P. Ferrari,
 M. Primicerio, G. Jona-Lasinio, S. Polidoro, J. Quastel   for many useful suggestions and comments. G. Carinci and C. Giardin\`a  acknowledge kind hospitality at the Universit\`a di l'Aquila.  A. De Masi thanks very warm hospitality at the Technion - Israel Institute of Technology.
 The research has been partially supported by PRIN 2009 (prot.
2009TA2595-002) and FIRB 2010 (grant n. RBFR10N90W).

\end{document}